\documentclass[10pt,fleqn]{amsart}

\usepackage{amsmath,amssymb,latexsym}
\usepackage[mathscr]{eucal}

\theoremstyle{plain}
\newtheorem{theorem}{Theorem}
\newtheorem{lemma}{Lemma}

\newcounter{rmkctr}
\newcommand{\rmk}{\refstepcounter{rmkctr}\smallskip\noindent{\bf\thermkctr. }}

\numberwithin{equation}{section}

\def\mydate{\number\year-\ifnum\month<10{0}\fi\number\month-\ifnum\day<10{0}\fi\number\day}

\newcommand{\dy}{\partial}

\newcommand{\ddt}[1]{\frac{\mathrm{d}{#1}}{\mathrm{d}{t}}}
\newcommand{\ddtau}[1]{\frac{\mathrm{d}{#1}}{\mathrm{d}{\tau}}}

\newcommand{\sfrac}[2]{{\textstyle\frac{#1}{#2}}}
\newcommand{\tssum}{{\textstyle\sum}}

\newcommand{\Zahl}{\mathbb{Z}}
\newcommand{\Real}{\mathbb{R}}
\newcommand{\Comp}{\mathbb{C}}

\newcommand{\ex}{\mathrm{e}}
\newcommand{\im}{\mathrm{i}}
\newcommand{\eps}{\varepsilon}
\newcommand{\vfi}{\varphi}

\newcommand{\eikx}{\ex^{\im\kb\cdot\xb}}

\newcommand{\dtau}{\>\mathrm{d}\tau}
\newcommand{\dz}{\>\mathrm{d}z}
\newcommand{\dx}{\>\mathrm{d}\boldsymbol{x}}

\newcommand{\lapl}[1]{\Delta_{#1}^{}}
\newcommand{\ilapl}[1]{\Delta_{#1}^{-1}}

\newcommand{\sump}{\mathop{\smash{\mathop{{\sum}'}}{\vphantom\sum}}}

\newcommand{\xb}{\boldsymbol{x}}
\newcommand{\ub}{\boldsymbol{u}}
\newcommand{\vb}{\boldsymbol{v}}
\newcommand{\ptil}{\delta p}

\newcommand{\ZL}{\Zahl_L}
\newcommand{\jb}{{\boldsymbol{j}}}

\newcommand{\kb}{{\boldsymbol{k}}}
\newcommand{\lb}{{\boldsymbol{l}}}

\newcommand{\Vb}{{\boldsymbol{V}}}

\newcommand{\gb}{\nabla}
\newcommand{\sgb}{\nabla^\perp}

\newcommand{\divv}{\gb\!\cdot\!\vb}

\newcommand{\rotv}{\sgb\!\cdot\!\vb}

\newcommand{\Pl}{{\sf P}^{\!{}^<}}

\newcommand{\Iw}{{\sf I}_\omega}

\newcommand{\Bw}{B_\omega}
\newcommand{\dst}{\dy_t^*}
\newcommand{\dstau}{\dy_\tau^*}

\newcommand{\Gy}{\mathcal{G}}

\newcommand{\Rm}{\mathcal{R}}
\newcommand{\Hm}{\mathcal{H}^\eps}

\newcommand{\Dom}{\mathscr{M}}

\newcommand{\fv}{f_{\vb}^{}}
\newcommand{\fr}{f_\rho^{}}

\newcommand{\cnst}[1]{c_{#1}^{}}
\newcommand{\cpoi}{\cnst{\textrm{p}}}
\newcommand{\cres}{\cnst{\textrm{nr}}}
\newcommand{\cU}{\cnst{U}}

\newcommand{\FD}{\mathsf{D}}

\DeclareMathOperator{\esssup}{ess\,sup}
\DeclareMathOperator{\sgn}{sgn}

\begin{document}

\title[Stability of the Slow Manifold]%
{Stability of the Slow Manifold\\ in the Primitive Equations}

\author[Temam]{R.~Temam}
\email{temam@indiana.edu}
\urladdr{http://mypage.iu.edu/\~{}temam}
\address[RT]{The Institute for Scientific Computing
   and Applied Mathematics\\
   Indiana University, Rawles Hall\\
   Bloomington, IN~47405--7106, United States}

\author[Wirosoetisno]{D.~Wirosoetisno}
\email{djoko.wirosoetisno@durham.ac.uk}
\urladdr{http://www.maths.dur.ac.uk/\~{}dma0dw}
\address[DW]{Department of Mathematical Sciences\\
   University of Durham\\
   Durham\ \ DH1~3LE, United Kingdom}

\thanks{This research was partially supported by the National Science
Foundation under grant NSF-DMS-0604235, by the Research Fund of Indiana
University, and by a grant from the Nuffield Foundation}

\keywords{Slow manifold, exponential asymptotics, primitive equations}
\subjclass[2000]{Primary: 35B40, 37L25, 76U05}


\begin{abstract}
We show that, under reasonably mild hypotheses, the solution of
the forced--dissipative rotating primitive equations of the ocean
loses most of its fast, inertia--gravity, component in the small
Rossby number limit as $t\to\infty$.
At leading order, the solution approaches what is known as ``geostrophic
balance'' even under ageostrophic, slowly time-dependent forcing.
Higher-order results can be obtained if one further assumes that
the forcing is time-independent and sufficiently smooth.
If the forcing lies in some Gevrey space, the solution will be exponentially
close to a finite-dimensional ``slow manifold'' after some time.
\end{abstract}

\maketitle


\section{Introduction}\label{s:intro}

One of the most basic models in geophysical fluid dynamics is the
primitive equations, understood here to be the hydrostatic approximation
to the rotating compressible Navier--Stokes equations, which is believed
to describe the large-scale dynamics of the atmosphere and the ocean to
a very good accuracy.
An important feature of such large-scale dynamics is that it largely
consists of slow motions in which the pressure gradient is nearly
balanced by the Coriolis force, a state known as {\em geostrophic balance\/}.
Various physical explanations have been given, some supported by numerical
simulations, to describe how this comes about, but to our knowledge
no rigorous mathematical proof has been proposed.
(For a review of the geophysical background, see, e.g., \cite{daley:ada}.)
One aim of this article is to prove that, in the limit of strong rotation
and stratification, the solution of the primitive equations will approach
geostrophic balance as $t\to\infty$, in the sense that the ageostrophic
energy will be of the order of the Rossby number.

As illustrated by the simple one-dimensional model \eqref{q:1dm},
here the basic mechanism for balance is the viscous damping of rapid
oscillations, leaving the slow dynamics mostly unchanged.
Separation of timescale, characterised by a small parameter $\eps$,
is therefore crucial for our result;
this is obtained by considering the limit of strong rotation {\em and\/}
stratification, or in other words, small Rossby number with Burger number
of order one.
We note that there are other physical mechanisms through which a balanced
state may be reached.
Working in an unbounded domain, an important example is the radiation of
inertia--gravity waves to infinity in what is known as the classical
geostrophic adjustment problem (see \cite[\S7.3]{gill:aod} and further
developments in \cite{reznik-al:01}).

Attempts to extend geostrophic balance to higher orders, and the closely
related problem of eliminating rapid oscillations in numerical solutions
(e.g., \cite{baer-tribbia:77,machenhauer:77,leith:80,vautard-legras:86}),
led naturally to the concept of {\em slow manifold\/} \cite{lorenz:86},
which has since become important in the study of rotating fluids (and
more generally of systems with multiple timescales).
We refer the reader to \cite{mackay:04} for a thorough review,
but for our purposes here, a slow manifold means a manifold in phase
space on which the normal velocity is small;
if the normal velocity is zero, we have an exact slow manifold.
In the geophysical literature, there have been many papers proposing
various formal asymptotic methods to construct slow manifolds
(e.g., \cite{warn-menard:86,wbsv:95}).
A number of numerical studies closely related to the stability of slow
manifolds have also been done (e.g., \cite{ford-mem-norton:00,polvani-al:94}).

It was realised early on \cite{lorenz:86,warn:97} that in general
no exact slow manifold exists and any construction is generally
asymptotic in nature.
For finite-dimensional systems, this can often be proved using
considerations of exponential asymptotics (see, e.g., \cite{kruskal-segur:91}).
More recently, it has been shown explicitly \cite{jv-yavneh:04} in an
infinite-dimensional rotating fluid model that exponentially weak
fast oscillations are generated spontaneously by vortical motion,
implying that slow manifolds could at best be exponentially accurate
(meaning the normal velocity on it be exponentially small).
Theorem~\ref{t:ho} shows, given the hypotheses, that an exponential accuracy
can indeed be achieved for the primitive equations, albeit with a weaker
dependence on $\eps$.

\medskip
From a more mathematical perspective, our exponentially slow manifold
(see Lemma~\ref{t:suba}), which is also presented in \cite{temam-dw:ebal} in
a slightly different form, is obtained using a technique adapted from
that first proposed in \cite{matthies:01}.
It involves truncating the PDE to a finite-dimensional system whose size
depends on $\eps$ and applying a classical estimate from perturbation
theory to the finite system.
By carefully balancing the truncation size and the estimates on
the finite system, one obtains a finite-dimensional exponentially accurate
slow manifold.
This estimate is local in time and only requires that the (instantaneous)
variables and the forcing be in some Sobolev space $H^s$;
it (although not the long-time asymptotic result below) can thus be
obtained for the inviscid equations as well.
If our solution is also Gevrey (which is true for the primitive equations
given Gevrey forcing),
the ignored high modes are exponentially small, so the ``total error''
(i.e.\ normal velocity on the slow manifold) is also exponentially small.

Gevrey regularity of the solution is therefore crucial in obtaining
exponential estimates.
As with the Navier--Stokes equations \cite{foias-temam:89}, in the
absence of boundaries and with Gevrey forcing, one can prove that
the strong solution of the primitive equations also has Gevrey
regularity \cite{petcu-dw:gev3}.
For the present article, we need uniform bounds on the norms, which
have been proved recently \cite{petcu:3dpe} following the global
regularity results of \cite{cao-titi:07,kobelkov:06,kobelkov:07}.
Since our result also assumes strong rotation, however, one could have
used an earlier work \cite{babin-al:00} which proved global regularity
under a sufficiently strong rotation and then used \cite{petcu-dw:gev3}
to obtain Gevrey regularity.

While our earlier paper \cite{temam-dw:ebal} is concerned with
a finite-time estimate on pointwise accuracy (``predictability''),
in this article our aim is to obtain long-time asymptotic estimates
(on ``balance'').
In this regard, the main problem for both the leading-order
(Theorem~\ref{t:o1}) and higher-order (Theorem~\ref{t:ho}) estimates
are the same: to bound the energy transfer, through the nonlinear term,
from the slow to fast modes at the same order as the fast modes themselves.
For this, one needs to handle not only {\em exact\/} fast--fast--slow
resonances, whose absence has long been known in the geophysical literature
(cf.\ e.g., \cite{bartello:95,embid-majda:96,lelong-riley:91,warn:86}
for discussions of related models), but also {\em near\/} resonances.
A key part in our approach is an estimate involving near resonances
in the primitive equations (cf.\ Lemma~\ref{t:nores}).
Another method based on algebraic geometry to handle related near resonances
can be found in \cite{babin-al:99}.

Taken together with \cite{temam-dw:ebal}, the results here may be regarded
as an extension of the single-frequency exponential estimates obtained in
\cite{matthies:01} to the ocean primitive equations, which have an infinite
number of frequencies.
Alternately, one may view Theorem~\ref{t:ho} as an extension to exponential
order of the leading-order results of \cite{babin-al:00} for a closely
related model.
Finally, our results here put a strong constraint on the nature of the
global attractor \cite{ju:07} in the strong rotation limit:
the attractor will have to lie within an exponentially thin
neighbourhood of the slow manifold.

\medskip
The rest of this article is arranged as follows.
We begin in the next section by describing the ocean primitive equations
(henceforth OPE) and recalling the known regularity results.
In Section~\ref{s:nm}, we write the OPE in terms of fast--slow variables
and in Fourier modes,
followed by computing explicitly the operator corresponding to the nonlinear
terms and describing its properties.
In Section~\ref{s:o1}, we state and prove our leading-order estimate, that
the solution of the OPE will be close to geostrophic balance as $t\to\infty$.
In the last section, we state and prove our exponential-order estimate.


\section{The Primitive Equations}\label{s:ope}

We start by recalling the basic settings of the ocean primitive equations
\cite{lions-temam-wang:92a}, and then recast the system in a form suitable
for our aim in this article.


\subsection{Setup}

We consider the primitive equations for the ocean,
scaled as in \cite{petcu-temam-dw:pe}
\begin{equation}\label{q:uvr}\begin{aligned}
   &\dy_t\vb + \frac1\eps \bigl[ \vb^\perp + \gb_2 p \bigr]
	+ \ub\cdot\gb \vb = \mu\Delta \vb + \fv,\\
   &\dy_t\rho - \frac1\eps u^3 + \ub\cdot\gb \rho = \mu\Delta \rho + \fr,\\
   &\gb\cdot\ub = \divv + \dy_z u^3 = 0,\\
   &\rho = -\dy_zp.
\end{aligned}\end{equation}
Here $\ub=(u^1,u^2,u^3)$ and $\vb=(u^1,u^2,0)$ are the three-
and two-dimensional fluid velocity, with $\vb^\perp:=(-u^2,u^1,0)$.
The variable $\rho$ can be interpreted in two ways:
One can take it to be the departure from a stably-stratified profile
(with the usual Boussinesq approximation), with the full density of
the fluid given by
\begin{equation}
   \rho_\textrm{full}(x,y,z,t) = \rho_0 - \eps^{-1} z\rho_1 + \rho(x,y,z,t),
\end{equation}
for some positive constants $\rho_0$ and $\rho_1$.
Alternately, one can think of it to be, e.g., salinity or temperature
that contributes linearly to the density.
The pressure $p$ is determined by the hydrostatic relation $\dy_zp=-\rho$
and the incompressibility condition $\gb\cdot\ub=0$,
and is not (directly) a function of $\rho$.
We write $\gb:=(\dy_x,\dy_y,\dy_z)$, $\gb_2:=(\dy_x,\dy_y,0)$,
$\Delta:=\dy_x^2+\dy_y^2+\dy_z^2$ and $\lapl2:=\dy_x^2+\dy_y^2$.
The parameter $\eps$ is related to the Rossby and Froude numbers;
in this paper we shall be concerned with the limit $\eps\to0$.
In general the viscosity coefficients for $\vb$ and $\rho$ are different;
we have set them both to $\mu$ for clarity of presentation (the general
case does not introduce any more essential difficulty).
The variables $(\vb,\rho)$ evidently depend on the parameters $\eps$ and $\mu$
as well as on $(\xb,t)$, but we shall not write this dependence explicitly.

We work in three spatial dimensions,
$\xb := (x,y,z) = (x^1,x^2,x^3) \in [0,L_1]\times[0,L_2]\times[-L_3/2,L_3/2]$\penalty0$=: \Dom$,
with periodic boundary conditions assumed;
we write $|\Dom|:=L_1L_2L_3$.
Moreover, following the practice in numerical simulations of stratified
turbulence (see, e.g., \cite{bartello:95}), we impose the following
symmetry on the dependent variables:
\begin{equation}\label{q:sym}\begin{aligned}
   &\vb(x,y,-z) = \vb(x,y,z),
	&\qquad
   &p(x,y,-z) = p(x,y,z),\\
   &u^3(x,y,-z) = -u^3(x,y,z),
	&\qquad
   &\rho(x,y,-z) = -\rho(x,y,z);
\end{aligned}\end{equation}
we say that $\vb$ and $p$ are {\em even} in $z$,
while $u^3$ and $\rho$ are {\em odd} in $z$.
For this symmetry to persist, $\fv$ must be even and $\fr$ odd in $z$.
Since $u^3$ and $\rho$ are also periodic in $z$,
we have $u^3(x,y,-L_3/2)=u^3(x,y,L_3/2)=0$
and $\rho(x,y,-L_3/2)=\rho(x,y,L_3/2)=0$;
similarly, $\dy_z u^1=0$, $\dy_z u^2=0$ and $\dy_z p=0$ on $z=0,\pm L_3/2$
if they are sufficiently smooth (as will be assumed below).
One may consider the symmetry conditions \eqref{q:sym} as a way to
impose the boundary conditions $u^3=0$, $\rho=0$, $\dy_z u^1=0$, $\dy_z u^2=0$
and $\dy_z p=0$ on $z=0$ and $z=L_3/2$ in the {\em effective domain}
$[0,L_1]\times[0,L_2]\times[0,L_3/2]$.
All variables and the forcing are assumed to have zero mean in $\Dom$;
the symmetry conditions above ensure that this also holds for their products
that appear below.
It can be verified that the symmetry \eqref{q:sym} is preserved by
the OPE \eqref{q:uvr}; that is, if it holds at $t=0$, it continues
to hold for $t>0$.


\subsection{Determining the pressure and vertical velocity}

Since $u^3=0$ at $z=0$, we can use (\ref{q:uvr}c) to write
\begin{equation}\label{q:u3}
   u^3(x,y,z) = -\int_0^z \divv(x,y,z') \>\mathrm{d}z'.
\end{equation}
Similarly, the pressure $p$ can be written in terms of the density $\rho$
as follows (cf.\ \cite{samelson-temam-wang2:03}).
Let $p(x,y,z)=\langle p(x,y)\rangle+\ptil(x,y,z)$ where $\langle\cdot\rangle$
denotes $z$-average and where
\begin{equation}\label{q:ptil}
   \ptil(x,y,z) = -\int_{z_0}^z \rho(x,y,z') \dz'
\end{equation}
with $z_0(x,y)$ chosen such that $\langle{\ptil}\rangle=0$;
this is most conveniently done using Fourier series (see below).
Using the fact that
\begin{equation}
   \int_{-L_3/2}^{L_3/2} \divv \dz
	= -\int_{-L_3/2}^{L_3/2} \dy_z u^3 \dz
	= u^3(\cdot,-L_3/2) - u^3(\cdot,L_3/2) = 0,
\end{equation}
and taking 2d divergence of the momentum equation (\ref{q:uvr}a),
we find
\begin{equation}
   \frac1\eps \bigl[ \gb\cdot\langle\vb^\perp\rangle + \lapl2 \langle p\rangle \bigr]
	+ \gb\cdot\langle{\ub\cdot\gb\vb}\rangle
	= \mu\Delta \gb\cdot\langle\vb\rangle + \gb\cdot\langle\fv\rangle.
\end{equation}
Here we have used the fact that $z$-integration commutes with horizontal
differential operators.
We can now solve for the average pressure $\langle p\rangle$,
\begin{equation}
   \langle p\rangle = \ilapl2\bigl[ -\gb\cdot\langle\vb^\perp\rangle + \eps \bigl(
	- \gb\cdot\langle\ub\cdot\gb\vb\rangle
	+ \mu\Delta\gb\cdot\langle\vb\rangle + \gb\cdot\langle\fv\rangle \bigr) \bigr]
\end{equation}
where $\ilapl2$ is uniquely defined to have zero $xy$-average.
With this, the momentum equation now reads
\begin{equation}\label{q:vb}\begin{aligned}
   \dy_t\vb + \frac1\eps\bigl[ \vb^\perp - \gb\ilapl2\gb\cdot\langle\vb^\perp\rangle
	&+ \gb_2\ptil \bigr] + \ub\cdot\gb\vb
	- \gb\ilapl2\gb\cdot\langle\ub\cdot\gb\vb\rangle\\
	&\hskip-10pt= \mu \Delta \bigl(\vb - \gb\ilapl2\gb\cdot\langle\vb\rangle\bigr)
	+ \fv - \gb\ilapl2\gb\cdot\langle\fv\rangle.
\end{aligned}\end{equation}


\subsection{Canonical form and regularity results}

Besides the usual $L^p(\Dom)$ and $H^s(\Dom)$, with $p\in[1,\infty]$
and $s\ge0$, we shall also need
the Gevrey space $G^\sigma(\Dom)$, defined as follows.
For $\sigma\ge0$, we say that $u\in G^\sigma(\Dom)$ if
\begin{equation}\label{q:Gevdef}
   |\ex^{\sigma(-\Delta)^{1/2}}u|_{L^2}^{} =: |u|_{G^\sigma}^{} < \infty.
\end{equation}
Let us denote our state variable $W=(\vb,\rho)^\mathrm{T}$.
We write $W\in L^p(\Dom)$ if $\vb\in L^p(\Dom)^2$, $\rho\in L^p(\Dom)$,
$(\vb,\rho)$ has zero average over $\Dom$ and $(\vb,\rho)$ satisfies
the symmetry \eqref{q:sym}, in the distribution sense as appropriate;
analogous notations are used for $W\in H^s(\Dom)$ and $W\in G^\sigma(\Dom)$,
and for the forcing $f$ (which has to preserve the symmetries of $W$).

With $u^3$ given by \eqref{q:u3} and $\ptil$ by \eqref{q:ptil},
we can write the OPE (\ref{q:uvr}b) and \eqref{q:vb} in the compact form
\begin{equation}\label{q:dW}
   \dy_t{W} + \frac1\eps LW + B(W,W) + AW = f.
\end{equation}
The operators $L$, $B$ and $A$ are defined by
\begin{equation}\begin{aligned}
   &LW = \bigl(\vb^\perp-\gb\ilapl2\gb\cdot\langle\vb^\perp\rangle+\gb_2\ptil,-u^3\bigr)^\textrm{T}\\
   &B(W,\hat W) = \bigl(\ub\cdot\gb\hat\vb
	- \gb\ilapl2\gb\cdot\langle\ub\cdot\gb\hat\vb\rangle,
	\ub\cdot\gb\hat\rho\bigr)^\textrm{T}\\
   &AW = -\bigl(\mu\Delta(\vb-\gb\ilapl2\gb\cdot\langle\vb\rangle),\mu\Delta\rho\bigr)^\textrm{T},
\end{aligned}\end{equation}
and the force $f$ is given by
\begin{equation}
   f = (\fv-\gb\ilapl2\gb\cdot\langle\fv\rangle,\fr)^\textrm{T}.
\end{equation}
The following properties are known (see, e.g., \cite{petcu-dw:gev3}).
The operator $L$ is antisymmetric: for any $W\in L^2(\Dom)$
\begin{equation}\label{q:Lasym}
   (LW,W)_{L^2} = 0;
\end{equation}
$B$ conserves energy: for any $W\in H^1(\Dom)$ and $\hat W\in H^1(\Dom)$,
\begin{equation}\label{q:Basym}
  (W,B(\hat W,W))_{L^2} = 0;
\end{equation}
and $A$ is coercive: for any $W\in H^2(\Dom)$,
\begin{equation}\label{q:Acoer}
   (AW,W)_{L^2}= \mu\,|\gb W|_{L^2}^2.
\end{equation}

We shall need the following regularity results for the OPE (here $K_s$
and $M_\sigma$ are continuous increasing functions of their arguments):

\newtheorem*{thmz}{Theorem 0}
\begin{thmz}\label{t:reg}
Let $W_0\in H^1$ and $f\in L^\infty(\Real_+;L^2)$.
Then for all $t\ge0$ there exists a solution $W(t)\in H^1$ of \eqref{q:dW}
with $W(0)=W_0$ and
\begin{equation}\label{q:WH1}
   |W(t)|_{H^1} \le K_0(|W_0|_{H^1},\|f\|_0^{})
\end{equation}
where, here and henceforth, $\|f\|_s^{}:=\esssup_{t\ge0}|f(t)|_{H^s}^{}$
for $s\ge0$.
Moreover, there exists a time $T_1(|W_0|_{H^1},\|f\|_0^{})$ such that
for $t\ge T_1$,
\begin{equation}\label{q:WH1u}
   |W(t)|_{H^1}^{} \le K_1(\|f\|_0^{}).
\end{equation}
Similarly, if $f\in L^\infty(\Real_+;H^{s-1})$, there exists a time
$T_s(|W_0|_{H^1},\|f\|_{s-1}^{})$ such that
\begin{equation}\label{q:WHsu}
   |W(t)|_{H^s}^{} \le K_s(\|f\|_{s-1}^{})
\end{equation}
for $t\ge T_s$.
Finally, fixing $\sigma>0$, if also
$\gb f\in L^\infty(\Real_+;G^\sigma)$, there exists a time
$T_\sigma(|W_0|_{H^1}^{},|\gb f|_{G^\sigma}^{})$ such that, for $t\ge T_\sigma$
\begin{equation}\label{q:WGsig}
   |\gb^2 W(t)|_{G^\sigma}^{} \le M_\sigma^{}(|\gb f|_{G^\sigma}).
\end{equation}
\end{thmz}

\noindent
The proof of \eqref{q:WH1}--\eqref{q:WH1u} can be found in \cite{ju:07};
the higher-order results \eqref{q:WHsu} can be found in \cite{petcu:3dpe}.
Both these works followed \cite{cao-titi:07} and \cite{kobelkov:06}.
The result \eqref{q:WGsig} follows from \cite{petcu-dw:gev3} and using
\eqref{q:WHsu} for $s=2$.

Since we are concerned with the limit of small $\eps$, however, one might also
be able to obtain \eqref{q:WH1} and \eqref{q:WHsu} following the method
used in \cite{babin-al:00} for the Boussinesq (non-hydrostatic) model.
One could then proceed to obtain \eqref{q:WGsig} as above.


\section{Normal Modes}\label{s:nm}

In this section, we decompose the solution $W$ into its slow and
fast components, expand them in Fourier modes, and state a lemma
that will be used in sections \ref{s:o1} and~\ref{s:ho} below.

\subsection{Fast and slow variables}

The Ertel potential vorticity
\begin{equation}
   q_E^{} = \rotv - \dy_z\rho
	+ \eps \bigl[(\dy_z\vb)\cdot\sgb\rho - \dy_z\rho\,(\rotv)\bigr],
\end{equation}
where $\sgb:=(-\dy_y,\dy_x,0)$, plays a central role in geophysical
fluid dynamics since it is a material invariant in the absence of
forcing and viscosity.
In this paper, however, it is easier to work with the {\em linearised\/}
potential vorticity (henceforth simply called {\em potential vorticity\/})
\begin{equation}
   q := \rotv - \dy_z\rho.
\end{equation}
From \eqref{q:uvr}, its evolution equation is
\begin{equation}
   \dy_t q + \sgb\cdot(\ub\cdot\gb\vb) - \dy_z(\ub\cdot\gb\rho) = \mu\Delta q
	+ f_q
\end{equation}
where $f_q:=\sgb\!\cdot\!\fv-\dy_z\fr$.
Let $\psi^0:=\Delta^{-1} q$, uniquely defined by requiring
that $\psi^0$ has zero integral over $\Dom$, and let
\begin{equation}\label{q:W0def}
    W^0 := \left(\begin{matrix} \vb^0\\ \rho^0 \end{matrix}\right)
	:= \left( \begin{matrix}\sgb\psi^0\\
		-\dy_z\psi^0 \end{matrix}\right).
\end{equation}
We note a mild abuse of notation on $\vb^0$ and $\sgb$:
$W^0=(-\dy_y\psi^0,\dy_x\psi^0,-\dy_z\psi^0)^\mathrm{T}$.

A little computation shows that $W^0$ lies in the kernel of the antisymmetric
operator $L$, that is, $LW^0=0$.
Conversely, if $LW=0$, then $W=(\sgb\Psi,-\dy_z\Psi)^\mathrm{T}$ for
some $\Psi$:
Since $u^3=0$, we have $\divv=0$, so
$\vb=\sgb\Psi+\Vb$ for some $\Psi(x,y,z)$ and $\Vb(z)$.
Now
\begin{equation}\label{q:kerL1}\begin{aligned}
   0 &= \vb^\perp - \gb_2\ilapl2\gb\!\cdot\!\langle v^\perp\rangle + \gb_2\delta p\\
     &= -\gb_2\Psi + \Vb^\perp + \gb_2\ilapl2\lapl2\langle\Psi\rangle
	+ \gb_2\delta p.
\end{aligned}\end{equation}
Since all other terms are horizontal gradients and $\Vb$ does
not depend on $(x,y)$, we must have $\Vb=0$.
Writing $\Psi(x,y,z)=\tilde\Psi(x,y,z)+\langle\Psi(x,y)\rangle$ where
$\tilde\Psi(x,y,z)$ has zero $z$-average, the terms that do not depend
on $z$ cancel and we are left with
\begin{equation}\label{q:kerL2}
   -\gb_2\tilde\Psi + \gb_2\delta p = 0.
\end{equation}
So $\delta p(x,y,z) = \tilde\Psi(x,y,z) + \Phi(z)$;
but since $\langle\delta p\rangle=0$, $\Phi=0$ and thus
$\rho=-\dy_z\Psi$ by \eqref{q:ptil}.
Therefore the null space of $L$ is completely characterised by \eqref{q:W0def},
\begin{equation}\label{q:kerL}
   \mathrm{ker}\,L = \{W^0:W^0=(\sgb\psi^0,-\dy_z\psi^0)^\mathrm{T}\}.
\end{equation}
With $\psi^0=\ilapl{}(\sgb\!\cdot\!\vb-\dy_z\rho)$ as above, this also
defines a projection $W\mapsto W^0$.
We call $W^0$ our {\em slow variable}.

Letting $B^0$ be the projection of $B$ to $\mathrm{ker}\,L$,
\begin{equation}
   B^0(W,\hat W) := \left( \begin{matrix}\sgb\Delta^{-1}\bigl[
	\sgb\cdot(\ub\cdot\gb\hat\vb) - \dy_z(\ub\cdot\gb\hat\rho)\bigr]\\
		-\dy_z\Delta^{-1}\vphantom{\Big|}\bigl[
	\sgb\cdot(\ub\cdot\gb\hat\vb) - \dy_z(\ub\cdot\gb\hat\rho)\bigr]
	\end{matrix}\right),
\end{equation}
we find that $W^0$ satisfies
\begin{equation}\label{q:dtW0}
   \dy_t{W^0} + B^0(W,W) + AW^0 = f^0
\end{equation}
where $f^0=(\sgb\Delta^{-1}f_q^{}, -\dy_z\Delta^{-1}f_q^{})^\mathrm{T}$ is the slow forcing.

Now let
\begin{equation}
   W^\eps = \left(\begin{matrix} \vb^\eps\\ \rho^\eps \end{matrix}\right)
	:= W - W^0
	= \left(\begin{matrix} \vb-\vb^0\\ \rho-\rho^0\end{matrix}\right).
\end{equation}
It will be seen below in Fourier representation that $W^\eps$ is a linear
combination of eigenfunctions of $L$ with imaginary eigenvalues whose moduli
are bounded from below;
we thus call $W^\eps$ our {\em fast variable}.
Since $\divv^0=0$, the vertical velocity $u^3$ is a purely fast variable.
In analogy with \eqref{q:dtW0}, we have
\begin{equation}\label{q:dtWeps}
   \dy_t{W^\eps} + \frac1\eps LW^\eps + B^\eps(W,W) + AW^\eps = f^\eps
\end{equation}
where $B^\eps(W,\hat W):=B(W,\hat W)-B^0(W,\hat W)$ and $f^\eps:=f-f^0$.

The fast variable has no potential vorticity, as can be seen by computing
$\sgb\cdot\vb^\eps-\dy_z\rho^\eps=q-\sgb\!\cdot\!\sgb\psi^0-\dy_{zz}\psi^0=0$.
Since the slow variable is completely determined by the potential vorticity,
this implies that the fast and slow variables are orthogonal in $L^2(\Dom)$,
\begin{equation}\label{q:Worth}\begin{aligned}
   (W^0,W^\eps)_{L^2}
	&= (\vb^0,\vb^\eps)_{L^2} + (\rho^0,\rho^\eps)_{L^2}\\
	&\hskip-3pt= (\sgb\psi^0,\vb^\eps)_{L^2} - (\dy_z\psi^0,\rho^\eps)_{L^2}
	= (\psi_0,-\sgb\cdot\vb^\eps+\dy_z\rho^\eps)_{L^2}
	= 0.
\end{aligned}\end{equation}

Of central interest in this paper is the ``fast energy''
\begin{equation}
   \sfrac12 |W^\eps|_{L^2}^2
	= \sfrac12\bigl(|\vb^\eps|_{L^2}^2+|\rho^\eps|_{L^2}^2\bigr).
\end{equation}
Its time derivative can be computed as follows.
Using \eqref{q:Worth}, we have after integrating by parts
\begin{equation}
   (W^\eps,\dy_tW)_{L^2} = (W^\eps,\dy_tW^0)_{L^2} + (W^\eps,\dy_tW^\eps)_{L^2}
	= \frac12\ddt{\:}|W^\eps|_{L^2}^2.
\end{equation}
Now \eqref{q:Basym} implies that
\begin{equation}
   (W^\eps,B(W,W))_{L^2} = (W^\eps,B(W,W^0+W^\eps))_{L^2} = (W^\eps,B(W,W^0))_{L^2}.
\end{equation}
Putting these together with \eqref{q:Lasym} and \eqref{q:Acoer}, we find
\begin{equation}\label{q:ddtweps}
   \frac12\ddt{}|W^\eps|_{L^2}^2 + \mu|\gb W^\eps|_{L^2}^2
	= -(W^\eps,B(W,W^0))_{L^2} + (W^\eps,f^\eps)_{L^2}.
\end{equation}


\subsection{Fourier expansion}

Thanks to the regularity results in Theorem~0, our solution $W(t)$ is
smooth and we can thus expand it in Fourier series,
\begin{equation}
   \vb(\xb,t) = \tssum_\kb^{}\, \vb_\kb(t)\, \eikx
   \qquad\textrm{and}\qquad
   \rho(\xb,t) = \tssum_\kb^{}\, \rho_\kb(t)\, \eikx.
\end{equation}
Here $\kb=(k_1,k_2,k_3)\in\ZL$ where
$\ZL=\Real^3/\Dom=\{(2\pi l_1/L_1,2\pi l_2/L_2, 2\pi l_3/L_3):(l_1,l_2,l_3)\in\Zahl^3\}$;
any wavevector $\kb$ is henceforth understood to live in $\ZL$.
We also denote $\kb':=(k_1,k_2,0)$ and write $\kb'\wedge\jb':=k_1j_2-k_2j_1$.
Since our variables have zero average over $\Dom$, $\vb_\kb=0$ when $\kb=0$;
moreover, since $\rho$ is odd in $z$, $\rho_\kb=0$ whenever $k_3=0$.
Thus $W_\kb:=(\vb_\kb,\rho_\kb)=0$ when $\kb=0$, which allows us to write
the $H^s$ norm simply as
\begin{equation}\label{q:Hsnorm}
   |W|_{H^s}^2 = \tssum_\kb^{}\,|\kb|^{2s}|W_\kb|^2
\end{equation}
and (see \eqref{q:Gevdef} for the definition of $G^\sigma$)
\begin{equation}\label{q:Gsignorm}
   |W|_{G^\sigma}^2 = \tssum_\kb^{}\,\ex^{2\sigma|\kb|}|W_\kb|^2\,.
\end{equation}

The antisymmetric operator $L$ is diagonal in Fourier space,
meaning that $L_{\kb\lb}=0$ when $\kb\ne\lb$;
we shall thus write $L_{\kb}:=L_{\kb\kb}$.
When $k_3\ne0$, we have
\begin{equation}
   L_\kb = \left( \begin{matrix} 0 &-1 &-k_1/k_3\\
			1 &0 &-k_2/k_3\\
			k_1/k_3 &k_2/k_3 &0 \end{matrix} \right).
\end{equation}
For $\kb'\ne0$, its eigenvalues are $\omega^0_\kb=0$ and
$\im\omega^\pm_\kb=\pm\im|\kb|/k_3$, where
$|\kb|:=\bigl(k_1^2+k_2^2+k_3^2)^{1/2}$, with eigenvectors
\begin{equation}\label{q:evectg}
   X^0_\kb = \frac1{|\kb|}\left(\begin{matrix} \phantom{-}k_2\\ -k_1\\ \phantom{-}k_3 \end{matrix}\right)
   \qquad\textrm{and}\qquad
   X^\pm_\kb = \frac1{\sqrt2|\kb'|\,|\kb|}\left(\begin{matrix}
	-k_2 k_3\pm\im k_1|\kb|\\ \phantom{-}k_1 k_3 \pm\im k_2|\kb|\\ |\kb'|^2
	\end{matrix}\right).
\end{equation}
When $\kb'=0$, we have $\omega_\kb^0=0$ and $\im\omega_\kb^\pm=\pm\im$
as eigenvalues with eigenvectors
\begin{equation}\label{q:evect0}
   X^0_\kb = \Biggl(\begin{matrix} \>0\>\\ 0\\ \sgn k_3 \end{matrix}\Biggr)
   \qquad\textrm{and}\qquad
   X^\pm_\kb = \frac{1}{\sqrt2}\Biggl(\begin{matrix} \>1\>\\ \mp\im\\ 0\end{matrix}\Biggr).
\end{equation}
For $\kb$ fixed, these eigenvectors are orthonormal under the inner product
$\cdot\;$ in $\Comp^3$.

When $k_3=0$, the fact that $\rho_\kb=0$ and $\kb\cdot\vb_\kb=0$
implies that the space is one-dimensional for each $\kb$ (in fact,
it is known that the vertically-averaged dynamics is that of the rotating
2d Navier--Stokes equations).
Since projecting to the $k_3=0$ subspace is equivalent to taking
vertical average, we compute
\begin{equation}
   \langle LW\rangle = (\langle\vb^\perp\rangle-\gb_2\ilapl2\gb\!\cdot\!\langle\vb^\perp\rangle,0)^\mathrm{T}
\end{equation}
where we have used $\langle u^3\rangle=0$ (since $u^3$ is odd) and
$\langle\delta p\rangle=0$ (by definition).
Reasoning as in \eqref{q:kerL1}--\eqref{q:kerL2} above, we find that
$\langle LW\rangle=0$, that is, the vertically-averaged ($k_3=0$) component is
completely slow.
In this case we can thus write
\begin{equation}\label{q:evect3}
   \omega_\kb^0 = 0
   \qquad\textrm{and}\qquad
   X^0_\kb = \frac1{|\kb'|} \left(\begin{matrix} \phantom{-}k_2\\ -k_1\\ \phantom{-}0 \end{matrix}\right),
\end{equation}
which can be included in the generic case $\kb'\ne0$ in computations.
Since the $k_3=0$ component is completely slow, $\langle W^\eps\rangle=0$,
there is no need to fix $X^\pm_\kb$.

We note that since $k_3\ne0$, $|\omega_\kb^\pm|\ge1$, viz.,
\begin{equation}\label{q:infw}
   \inf\, |\omega_\kb^\pm|^2 = \inf_{k_3\ne0}\, \biggl\{\frac{k_1^2+k_2^2+k_3^2}{k_3^2},\; 1\biggr\}
	= 1.
\end{equation}
In what follows, it is convenient to use $\{X^0_\kb,X^\pm_\kb\}$ as basis.

We can now write
\begin{equation}\label{q:W0Weps}\begin{aligned}
   &W^0(\xb,t) := \tssum_\kb\; w^0_\kb(t) X^0_\kb \eikx\\
   &W^\eps(\xb,t) := \tssum_\kb^s\; w^s_\kb(t) X^s_\kb
				\ex^{-\im\omega^s_\kb t/\eps}\eikx,
\end{aligned}\end{equation}
where $s\in\{-1,+1\}$, which we write as $\{-,+\}$ when it appears as a label.
The Fourier coefficients $w^0_\kb$ and $w^\pm_\kb$ are complex numbers
that depend on $t$ only, with $w^0_0=0$ and $w^\pm_{(k_1,k_2,0)}=0$.
With $\alpha\in\{-1,0,+1\}$, they can be computed using
\begin{equation}
   w^\alpha_\kb(t) = \frac1{|\Dom|} \int_\Dom W(\xb,t)\cdot X^\alpha_\kb
	\,\ex^{\im\omega^\alpha_\kb t/\eps-\im\kb\cdot\xb} \dx.
\end{equation}
The following relations hold:
\begin{equation}
   |W^0|_{L^2}^2 = \tssum_\kb\; |w^0_\kb|^2
   \qquad\textrm{and}\qquad
   |W^\eps|_{L^2}^2 = \tssum_\kb^s\; |w^s_\kb|^2.
\end{equation}
In addition, the fact that $(\vb^0,\rho^0)$ is real implies
\begin{equation}\label{q:w0}
   w^0_{-\kb} = -\overline{w^0_\kb}
   \qquad\textrm{and}\qquad
   w^0_{(k_1,k_2,-k_3)} = w^0_{(k_1,k_2,k_3)}
\end{equation}
where overbars denote complex conjugation.
Similarly, since $(\vb^\eps,\rho^\eps)$ is real,
\begin{equation}\label{q:weps1}
   w^\pm_{-\kb} = \overline{w^\pm_\kb}
   \qquad\textrm{and}\qquad
   w^\pm_{(k_1,k_2,-k_3)} = -w^\pm_{(k_1,k_2,k_3)}
\end{equation}
when $\kb'\ne0$ and, when $\kb'=0$,
\begin{equation}\label{q:weps2}
   w^\pm_{(0,0,-k_3)} = \overline{w^\mp_{(0,0,k_3)}}.
\end{equation}
We shall see below that, the linear oscillations having been factored out,
the variable $w^s_\kb$ is slow at leading order.
Similarly to $W$, we write the forcing $f$ as
\begin{equation}\begin{aligned}
   &f^0(\xb,t) := \tssum_\kb\; f^0_\kb(t) X^0_\kb \eikx\\
   &f^\eps(\xb,t) := \tssum_\kb^s\; f^s_\kb(t) X^s_\kb \eikx,
\end{aligned}\end{equation}
where, unlike in \eqref{q:W0Weps}, there is no factor of
$\ex^{-\im\omega_\kb^st/\eps}$ in the definition of $f^\eps$.
As noted above, $f$ must satisfy the same symmetries as $W$,
so the above properties of $w_\kb^\alpha$ also hold for $f_\kb^\alpha$;
we note in particular that $f_\kb^\pm=0$ when $k_3=0$.

For later convenience, we define the operator $\dy_t^*$ by
\begin{equation}
   \dy_t^* W := \ex^{-tL/\eps}\dy_t\,\ex^{tL/\eps}W.
\end{equation}
From \eqref{q:dW}, we find
\begin{equation}\label{q:dtSWeps}
   \dst W + B(W,W) + AW = f,
\end{equation}
which is $\dy_t W$ with the large antisymmetric term removed.

Now the nonlinear term on the rhs of \eqref{q:ddtweps} can be written as
\begin{equation}\begin{aligned}
   (W^\eps,B(W^0+W^\eps,W^0))_{L^2}
	&= (W^\eps,B(W^0,W^0))_{L^2} + (W^\eps,B(W^\eps,W^0))_{L^2}\\
	&= (W^\eps,B(W^0,W^0))_{L^2} - (W^0,B(W^\eps,W^\eps))_{L^2},
\end{aligned}\end{equation}
where the identity $(W^0,B(W^\eps,W^\eps))_{L^2}=-(W^\eps,B(W^\eps,W^0))_{L^2}$
had been obtained from \eqref{q:Basym}.

First, let

\vskip-15pt
\begin{equation}\begin{aligned}
   (W^\eps,B(W^0,W^0))_{L^2}
	&= |\Dom| \sum_{\jb\kb\lb}^s\, w^0_\jb w^0_\kb \overline{w^s_\lb}\,
		\im (X^0_\jb\cdot\kb')(X^0_\kb\cdot X^s_\lb)\,
			\delta_{\jb+\kb-\lb}\,\ex^{\im\omega^s_\lb t/\eps}\\
	&= \sum_{\jb\kb\lb}^s\, w^0_\jb w^0_\kb \overline{w^s_\lb}\,
		B_{\jb\kb\lb}^{00s} \ex^{\im\omega^s_\lb t/\eps}
\end{aligned}\end{equation}
where $\delta_{\jb+\kb-\lb}=1$ when $\jb+\kb=\lb$ and $0$ otherwise, and where
\begin{equation}\label{q:B00s}
   B_{\jb\kb\lb}^{00s} := \im\,|\Dom|\,\delta_{\jb+\kb-\lb}
	(X^0_\jb\cdot\kb')(X^0_\kb\cdot X^s_\lb).
\end{equation}
It is easy to verify from \eqref{q:B00s} that
$B_{\jb\kb\lb}^{00s}=0$ when $|\jb'|\,|\kb'|\,l_3=0$,
so we consider the other cases.
For the first factor, we have
\begin{equation}\label{q:B00sa}
   X^0_\jb\cdot\kb' = \frac{\kb'\wedge\jb'}{|\jb|}.
\end{equation}
For the second factor, we have
\begin{equation}\label{q:B00sb}\begin{aligned}
   &X^0_\kb\cdot X^s_\lb = \frac{k_2-\im s k_1}{\sqrt2\,|\kb|}
	&&\textrm{when } \lb'=0, \textrm{ and}\\
   &X^0_\kb\cdot X^s_\lb = \frac{k_3|\lb'|^2-(\kb'\cdot\lb')l_3-\im s(\lb'\wedge\kb')|\lb|}{\sqrt2\,|\kb|\,|\lb|\,|\lb'|}
	&&\textrm{when }\lb'\ne0.
\end{aligned}\end{equation}
From these, we have the bound
\begin{equation}\label{q:bdB00s}
   |B_{\jb\kb\lb}^{00s}| \le \frac{3\,|\Dom|}{\sqrt2}\frac{|\kb'|\,|\jb'|}{|\jb|}.
\end{equation}

Next, we consider
\begin{equation}\label{q:WBee0}\begin{aligned}
   (W^0,B(&W^\eps,W^\eps))_{L^2}\\
	&= |\Dom| \sum_{\jb\kb\lb}^{rs} w_\jb^r w_\kb^s \overline{w_\lb^0} \,\im\,({\sf V}X^r_\jb\cdot\kb)(X^s_\kb\cdot X^0_\lb) \,\delta_{\jb+\kb-\lb}\,\ex^{-\im(\omega_\jb^r+\omega_\kb^s)t/\eps}\\
	&= \sum_{\jb\kb\lb}^{rs}\, w_\jb^r w_\kb^s \overline{w_\lb^0}
		\, B_{\jb\kb\lb}^{rs0} \,\ex^{-\im(\omega_\jb^r+\omega_\kb^s)t/\eps}
\end{aligned}\end{equation}
where
\begin{equation}\label{q:Bee0}
   B_{\jb\kb\lb}^{rs0} := \im\,|\Dom|\,\delta_{\jb+\kb-\lb}\,({\sf V}X^r_\jb\cdot\kb)(X^s_\kb\cdot X^0_\lb)
\end{equation}
and where the operator ${\sf V}$, which produces an incompressible
velocity vector out of $X^r_\jb$, is defined by
\begin{equation}\begin{aligned}
   &{\sf V}X_\jb^r = X_\jb^r \hbox to120pt{}
	&&\textrm{when } j_3|\jb'|=0, \textrm{ and}\\
   &{\sf V}X_\jb^r = \frac1{\sqrt2\,|\jb|\,|\jb'|}\left(\begin{matrix}
	-j_2 j_3 +\im r j_1|\jb|\\
	\phantom{-}j_1 j_3 +\im r j_2|\jb|\\
	-\im r |\jb'|^2|\jb|/j_3  \end{matrix}\right)
	&&\textrm{when } j_3|\jb'|\ne 0.
\end{aligned}\end{equation}
Thus, we have ${\sf V}X^r_\jb\cdot\kb = 0$ when $j_3=0$,
\begin{equation}
   {\sf V}X^r_\jb\cdot\kb = \bigl(k_1-\im r k_2\bigr)/\sqrt2
\end{equation}
when $\jb'=0$, and
\begin{equation}
   {\sf V}X^r_\jb\cdot\kb = \frac{j_3(\jb'\wedge\kb') + \im r|\jb|(\jb'\cdot\kb') - \im r|\jb'|^2|\jb|\,k_3/j_3}{\sqrt2\,|\jb|\,|\jb'|}
\end{equation}
in the generic case $j_3|\jb'|\ne0$.
In all cases, we have the bound
\begin{equation}\label{q:bdvwr0}
   |{\sf V}X^r_\jb\cdot\kb|
	\le |\Dom|\,\bigl(\sqrt2\,|\kb'|+|\jb'|\,|k_3|/|j_3|\bigr).
\end{equation}
Next, $X_\kb^s\cdot X_\lb^0=0$ when $k_3=0$ or $\kb'=\lb'=0$, and
\begin{equation}\begin{aligned}
   &X_\kb^s\cdot X_\lb^0 = \frac{l_2+\im s l_1}{\sqrt2\,|\lb|}
	&&\textrm{when } \lb'\ne0 \textrm{ and } \kb'=0,\\
   &X_\kb^s\cdot X_\lb^0 = \sgn l_3\frac{|\kb'|}{\sqrt2\,|\kb|}
	&&\textrm{when } \lb'=0 \textrm{ and } \kb'\ne0,\\
   &X_\kb^s\cdot X_\lb^0 = \frac{-(\kb'\cdot\lb')k_3+\im s(\kb'\wedge\lb')|\kb|+|\kb'|^2l_3}{\sqrt2\,|\kb|\,|\kb'|\,|\lb|}
	&&\textrm{when } |\kb'|\,|\lb'|\,k_3\ne0.
\end{aligned}\end{equation}
These give us the bound
\begin{equation}
   |X_\kb^s\cdot X_\lb^0| \le \sqrt{5/2}
\end{equation}
in all cases and, together with \eqref{q:bdvwr0}, when $j_3\ne0$,
\begin{equation}\label{q:bdBrs0}
   |B_{\jb\kb\lb}^{rs0}|
	\le \sqrt5\,|\Dom|\,\bigl(|\kb'| + |\jb'|\,|k_3|/|j_3|\bigr).
\end{equation}
When $j_3k_3=0$ or $\lb=0$, we have $B_{\jb\kb\lb}^{rs0}=0$.

\subsection{Fast--Fast--Slow Resonances}
We first write \eqref{q:WBee0} as
\begin{equation}
   (W^0,B(W^\eps,W^\eps))_{L^2}
	= \frac12\sum_{\jb\kb\lb}^{rs}\, w_\jb^r w_\kb^s \overline{w_\lb^0}
	\,\bigl(B_{\jb\kb\lb}^{rs0} + B_{\kb\jb\lb}^{sr0}\bigr)
	\,\ex^{-\im(\omega_\jb^r+\omega_\kb^s)t/\eps}.
\end{equation}
It has long been known in the geophysical community that many rotating
fluid models ``have no fast--fast--slow resonances'' (see, e.g.,
\cite{warn:86} for the shallow-water equations and \cite{bartello:95}
for the Boussinesq equations).
In our notation, the absence of {\em exact\/} fast--fast--slow resonances
means that $B_{\jb\kb\lb}^{rs0}+B_{\kb\jb\lb}^{sr0}=0$ whenever
$\omega_\jb^r+\omega_\kb^s=0$;
the significance of this will be apparent below [see the development
following \eqref{q:Bst}].
For our purpose, however, we also need to consider {\em near\/} resonances,
i.e.\ those cases when $|\omega_\jb^r+\omega_\kb^s|$ is small but nonzero.
The following ``no-resonance'' lemma contains the estimate we need:

\begin{lemma}\label{t:nores}
For any $\jb$, $\kb$, $\lb\in\ZL$ with $\lb\ne0$,
\begin{equation}\label{q:nores}
   \bigl|B_{\jb\kb\lb}^{rs0}+B_{\kb\jb\lb}^{sr0}\bigr|
	\le \cres\,|\Dom|\,
		\Bigl(\frac{|\jb|\,|\kb|}{|\lb|} + |j_3| + |k_3|\Bigr)\,
		|\omega_\jb^r+\omega_\kb^s|
\end{equation}
where $\cres$ is an absolute constant.
\end{lemma}

\noindent
We note that $B_{\jb\kb\lb}^{rs0}=B_{\kb\jb\lb}^{sr0}=0$ when $\lb=0$
by \eqref{q:WBee0}, so this case is trivial.
We defer the proof to Appendix~\ref{s:nores}.



\section{Leading-Order Estimates}\label{s:o1}

In this section, we discuss the leading-order case of our general problem.
This is done separately due to its geophysical interest and since it
requires qualitatively weaker hypotheses.
As before, $W(t)=W^0(t)+W^\eps(t)$ is the solution of the OPE \eqref{q:dW}
with initial conditions $W(0)=W_0$, and
$K_{\rm g}(\cdot)$ is a continuous and increasing function of its argument.

\begin{theorem}\label{t:o1}
Suppose that the initial data $W_0\in H^1(\Dom)$ and that the forcing
$f\in L^\infty(\Real_+;H^2)\cap W^{1,\infty}(\Real_+;L^2)$, with
\begin{equation}
   \|f\|_{\rm g}^{} := \esssup_{t>0}\,\bigl(|f(t)|_{H^2}^{} + |\dy_tf(t)|_{L^2}^{}\bigr).
\end{equation}
Then there exist $T_{\rm g}=T_{\rm g}(|W_0|_{H^1}^{},\|f\|_{\rm g}^{},\eps)$ and $K_{\rm g}=K_{\rm g}(\|f\|_{\rm g}^{})$,
such that for $t\ge T_{\rm g}$,
\begin{equation}
   |W^\eps(t)|_{L^2}^{} \le \sqrt\eps\, K_{\rm g}(\|f\|_{\rm g}^{}).
\end{equation}
\end{theorem}

In geophysical parlance, our result states that, for given initial data
and forcing, the solution of the OPE will become geostrophically balanced
(in the sense that the ageostrophic component $W^\eps$ is of order
$\sqrt\eps$) after some time.
We note that the forcing may be time-dependent (although $\|f\|_{\rm g}^{}$ cannot
depend on $\eps$) and need not be geostrophic;
this will not be the case when we consider higher-order balance later.
Also, in contrast to the higher-order result in the next section,
no restriction on $\eps$ is necessary in this case.

\medskip
The linear mechanism of this ``geostrophic decay'' may be appreciated by
modelling \eqref{q:dtWeps}, without the non\-linear term, by the following ODE
\begin{equation}\label{q:1dm}
   \ddt{x} + \frac{\im}{\eps}\, x + \mu x = f
\end{equation}
where $\mu>0$ is a constant and $f=f(t)$ is given independently of $\eps$.
The skew-hermitian term $\im x/\eps$ causes oscillations of $x$ whose
frequency grows as $\eps\to0$.
In this limit, the forcing becomes less effective since $f$ varies slowly
by hypothesis while the damping remains unchanged, so $x$ will eventually
decay to the order of the ``net forcing'' $\sqrt\eps f$.
More concretely, let $z(t)=\ex^{\im t/\eps}x(t)$ and write \eqref{q:1dm} as
\begin{equation}\label{q:dzdt}
   \ddt{\;}\bigl(\ex^{\mu t/2}z\bigr)
	+ \frac{\mu}{2}\ex^{\mu t/2}z = \ex^{\mu t/2-\im t/\eps}f,
\end{equation}
from which it follows that
\begin{equation}
   \ddt{\;}\bigl(\ex^{\mu t/2}|z|^2\bigr)
	+ \mu\ex^{\mu t/2}|z|^2
	= 2\ex^{\mu t/2}\mathrm{Re}\,\bigl(\ex^{-\im t/\eps}\bar z f\bigr).
\end{equation}
Integrating, we find
\begin{equation}\begin{aligned}
   \ex^{\mu t/2}|z(t)| &- |z(0)|^2 + \mu\int_0^t \ex^{\mu\tau/2} |z(\tau)|^2\dtau
	= 2 \int_0^t \ex^{\mu\tau/2} \mathrm{Re}\bigl(\ex^{-\im\tau/\eps}\bar z f\bigr)\dtau\\
	&= 2\eps \bigl[\ex^{\mu\tau/2}\mathrm{Re}\bigl(\im\ex^{-\im t/\eps}\bar z f\bigr)\bigr]_0^t
	- 2\eps \int_0^t \mathrm{Re}\bigl[\im\ex^{-\im\tau/\eps}\dy_\tau(\ex^{\mu\tau/2}\bar z f)\bigr] \dtau,
\end{aligned}\end{equation}
where the second equality is obtained by integration by parts.
Since $\dy_tf$ is bounded independently of $\eps$, the integral can be
bounded using \eqref{q:dzdt} and the integral on the left-hand side.
This leaves us with
\begin{equation}
   |z(t)|^2 \le \ex^{-\mu t/2}\,\cnst1(|f|)\,|z(0)|^2
	+ \frac{\eps}{\mu}\,(1-\ex^{-\mu t/2})\,K(|f|,|\dy_tf|,\mu).
\end{equation}
Most of the work in the proof below is devoted to handling the nonlinear term,
where particular properties of the OPE come into play.
A PDE application of this principle can be found in \cite{schochet:94}.


\subsection{Proof of Theorem~\ref{t:o1}}

In this proof, we omit the subscript in the inner product
$(\cdot,\cdot)_{L^2}$ when the meaning is unambiguous;
similarly, $|\cdot|\equiv|\cdot|_{L^2}^{}$.
We start by writing \eqref{q:ddtweps} as
\begin{equation}\label{q:dt0Weps}\begin{aligned}
   \ddt{} |W^\eps|^2 &+ 2\mu |\gb W^\eps|^2\\ 
	&= -2(W^\eps,B(W^0,W^0)) - 2(W^\eps,B(W^\eps,W^0)) + 2(W^\eps,f^\eps)\\
	&= -2(W^\eps,B(W^0,W^0)) + 2(W^0,B(W^\eps,W^\eps)) + 2(W^\eps,f^\eps).
\end{aligned}\end{equation}
Using the Poincar{\'e} inequality, $|W^\eps|^2\le\cpoi|\gb W^\eps|^2$,
and multiplying the left-hand side by $2\ex^{\nu t}$ where $\nu:=\mu\cpoi$,
we have
\begin{equation}
   \ddt{}\bigl( \ex^{\nu t} |W^\eps|^2 \bigr) + \mu\ex^{\nu t} |\gb W^\eps|^2
	\le \ex^{\nu t}\Bigl(\ddt{} |W^\eps|^2 + \mu |\gb W^\eps|^2 + \mu |\gb W^\eps|^2\Bigr).
\end{equation}
With this, \eqref{q:dt0Weps} becomes
\begin{equation}\label{q:ddtwfour}\begin{aligned}
   &\ddt{\;} \bigl(\ex^{\nu t}\, |W^\eps|^2\bigr)
		+ \mu \ex^{\nu t}|\gb W^\eps|^2\\
	&\hbox to24pt{}\le 2\,\ex^{\nu t}\, (W^\eps,f^\eps) - 2\,\ex^{\nu t}\, (W^\eps,B(W^0,W^0)) + 2\,\ex^{\nu t}\, (W^0,B(W^\eps,W^\eps)).
\end{aligned}\end{equation}

We now integrate this inequality from $0$ to $t$.
On the left-hand side we have
\begin{equation}\label{q:lhs1}\begin{aligned}
   \int_0^t \Bigl\{ \ddtau{\;} \bigl(\ex^{\nu\tau} |W^\eps|^2\bigr)
		&+ \mu \ex^{\nu\tau}|\gb W^\eps|^2 \Bigr\} \dtau\\
   &= \ex^{\nu t}|W^\eps(t)|^2 - |W^\eps(0)|^2
	+ \mu \int_0^t \ex^{\nu\tau}|\gb W^\eps|^2 \dtau.
\end{aligned}\end{equation}

Using the expansion \eqref{q:W0Weps} of $W^\eps$,
we integrate the right-hand side by parts to bring out a factor of $\eps$;
that is, we integrate the rapidly oscillating exponential
$\ex^{\im\omega_\kb^st/\eps}$ and leave everything else.
For the force term, we have
\begin{equation}\begin{aligned}
   \int_0^t \ex^{\nu\tau} (W^\eps,f^\eps) \dtau
	&= |\Dom| \sum_\kb^s\,\int_0^t \ex^{\nu\tau+\im\omega_\kb^s\tau/\eps}
		\overline{w_\kb^s} f_\kb^s \dtau\\
	&= \eps\, |\Dom|\sump_\kb^s\, \frac{1}{\im\omega_\kb^s}\bigl[\overline{w_\kb^s(t)} f_\kb^s(t) \ex^{\nu t+\im\omega_\kb^st/\eps} - \overline{w_\kb^s(0)} f_\kb^s(0)\bigr]\\
	&\qquad- \eps\, |\Dom|\int_0^t \sump_\kb^s\, \frac{\ex^{\im\omega_\kb^s\tau/\eps}}{\im\omega_\kb^s}
		\ddtau{\;}\bigl( \overline{w_\kb^s} f_\kb^s \ex^{\nu\tau}\bigr) \dtau.
\end{aligned}\end{equation}
Here the prime on $\sump$ indicates that terms for which $\omega_\kb^s=0$
are omitted since then $w_\kb^s=0$.
Introducing the integration operator $\Iw$ defined by
\begin{equation}\label{q:Iwdef}
   \Iw W^\eps(\xb,t) := \sump_\kb^s\; \frac\im{\omega^s_\kb}\, w^s_\kb(t) X^s_\kb
				\ex^{-\im\omega^s_\kb t/\eps}\eikx\,,
\end{equation}
which is well-defined since $|\omega_\kb^s|\ge1$, we can write this as
\begin{equation}\label{q:Wef}\begin{aligned}
   \int_0^t &\ex^{\nu\tau} (W^\eps,f^\eps) \dtau
	= \eps\, \ex^{\nu t}(\Iw W^\eps(t),f^\eps(t))
	- \eps\,(\Iw W^\eps(0),f^\eps(0))\\
	&- \eps \int_0^t \ex^{\nu\tau}\bigl\{ \nu(\Iw W^\eps,f^\eps)
		+ (\Iw\dstau W^\eps,f^\eps)
		+ (\Iw W^\eps,\dy_\tau f^\eps) \bigr\} \dtau.
\end{aligned}\end{equation}
Similarly, integrating the next term by parts we find
\begin{equation}\label{q:WzBze}\begin{aligned}
   &\int_0^t \ex^{\nu\tau} (W^\eps,B(W^0,W^0)) \dtau\\
	&\quad= \eps\, \ex^{\nu t}(\Iw W^\eps,B(W^0,W^0))(t)
	- \eps\, (\Iw W^\eps,B(W^0,W^0))(0)\\
	&\qquad- \eps \int_0^\tau \ex^{\nu\tau} \bigl\{ \nu\, (\Iw W^\eps,B(W^0,W^0))
		+ (\Iw\dstau W^\eps,B(W^0,W^0))\\
		&\hbox to75pt{}+ (\Iw W^\eps,B(\dy_\tau W^0,W^0))
		+ (\Iw W^\eps,B(W^0,\dy_\tau W^0)) \bigr\}\dtau.
\end{aligned}\end{equation}

Next, we consider
\begin{equation}\label{q:Bst}\begin{aligned}
   \int_0^t &\ex^{\nu\tau}\,(W^0,B(W^\eps,W^\eps))\,\dtau\\
	&= \int_0^t \frac12\sum_{\jb\kb\lb}^{rs}
		\ex^{-\im(\omega_\jb^r+\omega_\kb^s)\tau/\eps}
		\bigl(B_{\jb\kb\lb}^{rs0} + B_{\kb\jb\lb}^{sr0}\bigr)
		w_\jb^r w_\kb^s \overline{w_\lb^0}\, \ex^{\nu\tau} \dtau\\
	&= \frac{\eps\im}2\,\sump_{\jb\kb\lb}^{rs}
		\frac{B_{\jb\kb\lb}^{rs0}+B_{\kb\jb\lb}^{sr0}}{\omega_\jb^r+\omega_\kb^s}\,
		[w_\jb^r(t)w_\kb^s(t)\overline{w_\lb^0(t)}\ex^{\nu t-\im(\omega_\jb^r+\omega_\kb^s)t/\eps}\\
	&\hbox to210pt{}- w_\jb^r(0)w_\kb^s(0)\overline{w_\lb^0(0)}]\\
	&\qquad {}- \frac{\eps\im}2\int_0^t \sump_{\jb\kb\lb}^{rs}
		\frac{B_{\jb\kb\lb}^{rs0}+B_{\kb\jb\lb}^{sr0}}{\omega_\jb^r+\omega_\kb^s}\,
		\ex^{-\im(\omega_\jb^r+\omega_\kb^s)\tau/\eps}
		\ddtau{}\bigl[ w^r_\jb w^s_\kb \overline{w^0_\lb} \ex^{\nu \tau}\bigr]\;\mathrm{d}\tau.
\end{aligned}\end{equation}
Here the prime on $\sum'$ indicates that exactly resonant terms,
for which $\omega_\jb^r+\omega_\kb^s=0$ and
$B_{\jb\kb\lb}^{rs0}+B_{\kb\jb\lb}^{sr0}=0$, are excluded.
Using the bilinear operator $\Bw$,
defined for any $W^\eps$, $\hat W^\eps$ and $\tilde W^0$ by
\begin{equation}\label{q:Bwsdef}
   (\tilde W^0,\Bw(W^\eps,\hat W^\eps))
	:= \frac{\im}2 \sump_{\jb\kb\lb}^{rs}
	\frac{B_{\jb\kb\lb}^{rs0}+B_{\kb\jb\lb}^{sr0}}{\omega_\jb^r+\omega_\lb^s}
		w_\jb^r \hat w_\kb^s \overline{\tilde w\vphantom{w}_\lb^0}
		\ex^{-\im(\omega_\jb^r+\omega_\kb^s)t/\eps},
\end{equation}
we can write \eqref{q:Bst} in the more compact form
\begin{equation}\label{q:Bws}\begin{aligned}
   \!\int_0^t &\ex^{\nu\tau}\, (W^0,B(W^\eps,W^\eps)) \dtau\\
	&= \eps\,\ex^{\nu t}\,(W^0,\Bw(W^\eps,W^\eps))(t)
	- \eps\,(W^0,\Bw(W^\eps,W^\eps))(0)\\
	&\quad {}- \eps\int_0^t \ex^{\nu\tau}\,\bigl\{
		\nu\,(W^0,\Bw(W^\eps,W^\eps))
	+ (\dy_\tau W^0,\Bw(W^\eps,W^\eps))\\
	&\hskip156pt {}+ (W^0,\dstau\Bw(W^\eps,W^\eps)) \}\dtau.
\end{aligned}\end{equation}

Putting these together, \eqref{q:ddtwfour} integrates to
\begin{equation}\label{q:Wt1}\begin{aligned}
   \ex^{\nu t}|W^\eps(t)|^2 &- |W^\eps(0)|^2 + \mu \int_0^t \ex^{\nu\tau}|\gb W^\eps|^2 \dtau\\
	&\le 2\eps\,\ex^{\nu t}\, \bigl(\Iw W^\eps,f^\eps\bigr)(t)
		- 2\eps\, \bigl(\Iw W^\eps,f^\eps\bigr)(0)\\
	&\quad- 2\eps\, \ex^{\nu t}\, \bigl(\Iw W^\eps,B(W^0,W^0)\bigr)(t)
	+ 2\eps\, \bigl(\Iw W^\eps,B(W^0,W^0)\bigr)(0)\\
	&\quad{}+ 2\eps\,\ex^{\nu t}\bigl(W^0,\Bw(W^\eps,W^\eps)\bigr)(t)
		- 2\eps\,\bigl(W^0,\Bw(W^\eps,W^\eps)\bigr)(0)\\
	&\quad+ 2\eps \int_0^t \ex^{\nu\tau} \bigl\{
		I_0(\tau) - I_1(\tau) + I_2(\tau)\bigr\} \dtau.
\end{aligned}\end{equation}
Here the integrands are
\begin{equation}
   I_0 := \nu (\Iw W^\eps,f^\eps) + (\Iw W^\eps,\dy_\tau f^\eps) + (\Iw\dstau W^\eps,f^\eps),
\end{equation}
\begin{equation}\begin{aligned}
   I_1 &:= \nu\, (\Iw W^\eps, B(W^0,W^0))
        + (\Iw\dstau W^\eps,B(W^0,W^0))\\
	&\qquad {}+ (\Iw W^\eps, B(\dy_\tau W^0,W^0))
	+ (\Iw W^\eps,B(W^0,\dy_\tau W^0)),
\end{aligned}\end{equation}
and
\begin{equation}\begin{aligned}
   I_2 := \nu\,(W^0,\Bw(W^\eps,W^\eps))
	&+ (\dy_\tau W^0,\Bw(W^\eps,W^\eps))\\
	&+ (W^0,\dstau\Bw(W^\eps,W^\eps)).
\end{aligned}\end{equation}

We now bound the right-hand side of \eqref{q:Wt1}.
On the second line, we have
\begin{equation}\begin{aligned}
   \bigl|\ex^{\nu t}\,(\Iw W^\eps(t),f^\eps(t))&-(\Iw W^\eps(0),f^\eps(0))\bigr|\\
	&\le \ex^{\nu t}\,|W^\eps(t)|\,|f^\eps(t)| + |W^\eps(0)|\,|f^\eps(0)|,
\end{aligned}\end{equation}
where we have used the fact that, thanks to \eqref{q:infw},
\begin{equation}
   |\gb^\alpha \Iw W^\eps| \le |\gb^\alpha W^\eps|,
   \qquad\textrm{for }\alpha=0,1,2,\cdots.
\end{equation}
To bound the next line, we use the estimate
\begin{equation}\label{q:B0ee}\begin{aligned}
   |(\tilde W,B(W^0,\hat W))|
	&\le C\, |\tilde W|_{L^6}\,|W^0|_{L^3}\,|\gb\hat W|_{L^2}\\
	&\le C\, |\gb\tilde W|\,|W^0|^{1/2}\,|\gb W^0|^{1/2}\,|\gb\hat W|
\end{aligned}\end{equation}
(note that the first argument of $B$ is $W^0$) to obtain
\begin{equation}\begin{aligned}
   \bigl|\ex^{\nu t}\,(\Iw W^\eps,B(W^0,\,&W^0))(t)-(\Iw W^\eps,B(W^0,W^0)(0)\bigr|\\
	&\le \ex^{\nu t}\,|\gb W^\eps(t)|\,|W^0(t)|^{1/2}|\gb W^0(t)|^{3/2}\\
	&\qquad+ |\gb W^\eps(0)|\,|W^0(0)|^{1/2}|\gb W^0(0)|^{3/2}.
\end{aligned}\end{equation}
In \eqref{q:B0ee} and in the rest of this proof, $C$ and $c$ denote generic
constants which may not be the same each time the symbol is used; such
constants may depend on $\Dom$ but not on any other parameter.
Numbered constants may also depend on $\mu$.

We now derive a bound involving $\Bw$.
Since $B_{\jb\kb\lb}^{rs0} + B_{\kb\jb\lb}^{sr0}=0$ in the case of exact
resonance, we assume that $\omega_\jb^r+\omega_\kb^s\ne0$.
Then \eqref{q:nores} implies
\begin{equation}
   \frac{\bigl|B_{\jb\kb\lb}^{rs0} + B_{\kb\jb\lb}^{sr0}\bigr|}{|\omega_\jb^r+\omega_\kb^s|}
	\le C\,|\jb|\,|\kb|.
\end{equation}
With this, we have for any $W^\eps$, $\hat W^\eps$ and $\tilde W^0$,
\begin{equation}\label{q:bd3Bws}\begin{aligned}
   \!\bigl|(\tilde W^0,\Bw(W^\eps,\hat W^\eps))\bigr|
	&\le \frac12\sump_{\jb\kb\lb}^{rs}
	\biggl|\frac{B_{\jb\kb\lb}^{rs0}+B_{\kb\jb\lb}^{sr0}}{\omega_\jb^r+\omega_\kb^s}\biggr|\,
 		|w^r_\jb|\,|\hat w^s_\kb|\,|\tilde w_\lb^0|\\
	&\le C\,\sum_{\jb+\kb=\lb}^{rs}\,
		|\jb|\,|\kb|\,|w^r_\jb|\,|\hat w^s_\kb|\,|\tilde w_\lb^0|\\
	&\le \int_\Dom \theta(\xb)\,\xi(\xb)\,\zeta(\xb)\;\mathrm{d}\xb^3\\
	&\le C\,|\gb W^\eps|_{L^p}|\gb\hat W^\eps|_{L^q}|\tilde W^0|_{L^m}\,,
\end{aligned}\end{equation}
with $1/p+1/q+1/m=1$ and where on the penultimate line
\begin{equation}
   \theta(\xb) := \tssum_\jb^r\,|\jb|\,|w_\jb^r|\,\ex^{\im\jb\cdot\xb},
   \>
   \xi(\xb) := \tssum_\kb^s\,|\kb|\,|\hat w_\kb^s|\,\ex^{\im\kb\cdot\xb}
   \textrm{ and }
   \zeta(\xb) := \tssum_\lb\,|\tilde w_\lb^0|\,\ex^{\im\lb\cdot\xb}.
\end{equation}

Using \eqref{q:bd3Bws} with $p=q=2$ and $m=\infty$, plus the embedding
$H^2\subset\subset L^\infty$, we have the bound
\begin{equation}\begin{aligned}
   \bigl|\ex^{\nu t}\bigl(W^0,\,&\Bw(W^\eps,W^\eps)\bigr)(t) - \bigl(W^0,\Bw(W^\eps,W^\eps)\bigr)(0)\bigr|\\
	&\le C\,\bigl(\ex^{\nu t}|\gb W^\eps(t)|^2|\gb^2W^0(t)|
		+ |\gb W^\eps(0)|^2|\gb^2W^0(0)|\bigr).
\end{aligned}\end{equation}

To bound the integrand in \eqref{q:Wt1}, we need estimates on $\dy_tW^0$ and
$\dst W^\eps$ in addition to those already obtained.
Using the bound
\begin{equation}
   |B^0(W,W)|_{L^2}^{} \le C\,|\gb W|_{L^4}^2
	\le C\,|\gb W|_{H^{3/4}}^2
	\le C\,|\gb^2 W|^{3/2}|\gb W|^{1/2},
\end{equation}
we find from \eqref{q:dtW0}
\begin{equation}\label{q:bd0dtW0}\begin{aligned}
   |\dy_t W^0|_{L^2}^{}
	&\le C\,|\gb W|_{H^{3/4}}^2 + \mu\,|\gb^2 W^0| + |f|\\
	&\le C\,|\gb^2 W|^{3/2}|\gb W|^{1/2} + \mu\, |\gb^2 W^0| + |f|.
\end{aligned}\end{equation}
Similarly, we find from \eqref{q:dtSWeps}
\begin{equation}\label{q:bd0dtWe}\begin{aligned}
   |\dst W^\eps|_{L^2}^{}
	&\le C\,|\gb W|_{H^{3/4}}^2 + \mu\,|\gb^2 W^\eps| + |f|\\
	&\le C\,|\gb^2 W|^{3/2}|\gb W|^{1/2} + \mu\, |\gb^2 W^\eps| + |f|.
\end{aligned}\end{equation}
Now using the bound
\begin{equation}
   |\gb B(W,W)|_{L^2}^{} \le C\,|\gb^2 W|_{L^{12/5}}^{}|\gb W|_{L^{12}}^{}
	\le C\,|\gb^2 W|_{H^{1/4}}^{}
\end{equation}
we find
\begin{equation}\label{q:bd1dtW}\begin{aligned}
   |\gb\dy_t W^0|_{L^2}^{}
	&\le C\,|\gb^2 W|_{H^{1/4}}^2 + \mu\,|\gb^3 W^0| + |\gb f^0|\\
	&\le C\,|\gb^3 W|^{1/2}|\gb^2 W|^{3/2} + \mu\,|\gb^3 W^0| + |\gb f^0|,\\
   |\gb\dst W^\eps|_{L^2}^{}
	&\le C\,|\gb^2 W|_{H^{1/4}}^2 + \mu\,|\gb^3 W^\eps| + |\gb f^\eps|\\
	&\le C\,|\gb^3 W|^{1/2}|\gb^2 W|^{3/2} + \mu\,|\gb^3 W^\eps| + |\gb f^\eps|.\\
\end{aligned}\end{equation}

The bound for $I_0$ follows by using \eqref{q:bd0dtWe},
\begin{equation}\label{q:bdI0}\begin{aligned}
   |I_0|_{L^2}^{} &\le C\,\bigl(|\gb W|_{H^{3/4}}^2
	+ (\mu+c)\,|\gb^2W^\eps| + |f^\eps| \bigr)\,(|f^\eps|+|\dy_t f^\eps|)\\
	&\le C\,\bigl(|\gb W|_{H^{3/4}}^2
		+ (\mu+c)\,|\gb^2W| + \|f\|_{\rm g}^{} \bigr)\,\|f\|_{\rm g}^{},
\end{aligned}\end{equation}
where we have used the fact that
$|\gb^\alpha W^\eps|^2 \le |\gb^\alpha W^\eps|^2 + |\gb^\alpha W^0|^2 = |\gb^\alpha W|^2$.
Next, using \eqref{q:bd3Bws} we bound $I_2$ as
\begin{equation}\begin{aligned}
   |I_2|_{L^2}^{}
	&\le \mu c\,|W^0|_{L^\infty}^{}|\gb W^\eps|^2
	+ c\,|\dy_\tau W^0|\,|\gb W^\eps|_{L^4}^2\\
	&\hskip93pt {}+ c\,|W^0|_{L^\infty}^{}|\gb W^\eps|\,|\gb\dst W^\eps|\\
	&\le \mu c\,|\gb^2W^0|\,|\gb W^\eps|^2
	+ c\,\bigl(|\gb W|_{H^{3/4}}^2 + \mu\,|\gb^2 W^0| + |f^0|\bigr)
		|\gb W^\eps|_{H^{3/4}}^2\\
	&\hskip20pt {}+ c\,|\gb^2 W^0|\,|\gb W^\eps|\,
		\bigl(|\gb^2W|_{H^{1/4}}^2 + \mu\,|\gb^3W^\eps| + |\gb f^\eps|\bigr)\\
	&\le c\,|\gb W|\,|\gb^2W|\,|\gb^2W|_{H^{1/4}}^2
	+ \mu c\,|\gb^3W|\,|\gb^2W|\,|\gb W|\\
	&\hskip20pt {}+ |\gb^2W|^{3/2}|\gb W|^{1/2}\,\|f\|_{\rm g}^{}
\end{aligned}\end{equation}
where interpolation inequalities have been used for the last step.
The bound for $I_1$ is majorised by that for $I_2$.


Putting everything together, we have from \eqref{q:Wt1}
\begin{equation}\label{q:Wt2}\begin{aligned}
   \ex^{\nu t}\,&|W^\eps(t)|^2 - |W^\eps(0)|^2\\
	&\le \eps\,c_2\,\ex^{\nu t}\,|\gb^2 W(t)|\bigl(|\gb W(t)|^2 + \|f\|_{\rm g}^{}\bigr)
		+ \eps\,c_2\,|\gb^2 W_0|\bigl(|\gb W_0|^2 + \|f\|_{\rm g}^{}\bigr)\\
	&\quad {}+ \eps\,c_3 \int_0^t \ex^{\nu\tau} \bigl\{
		|W|_{H^1}^{}|W|_{H^2}^{}|W|_{H^{9/4}}^2 + \mu\,|W|_{H^3}^{}|W|_{H^2}^{}|W|_{H^1}^{}\\
		&\hskip80pt {}+ \bigl(|W|_{H^2}^{3/2}|W|_{H^1}^{1/2} + (\mu+c)\,|W|_{H^2}^{} + \|f\|_{\rm g}^{}\bigr)\,\|f\|_{\rm g}^{} \bigr\} \dtau.
\end{aligned}\end{equation}
Now by \eqref{q:WH1u} and \eqref{q:WHsu}, we can find $K_*(\|f\|_{\rm g}^{})$
and $T_*(|\gb W_0|,\|f\|_{\rm g}^{})$ such that, for $t\ge T_*$,
\begin{equation}
   c\,|\gb^sW(t)|^2 + (\mu+c')\, |\gb^s W(t)| + \|f\|_{\rm g}^{} \le K_*
\end{equation}
for $s\in\{0,1,2,3\}$.
Let $t':=t-T_*$ and relabel $t$ in \eqref{q:Wt2} as $t'$.
We can then bound the integral in \eqref{q:Wt2} as
\begin{equation}
   \int_0^{t'} \ex^{\nu\tau}\bigl\{\cdots\}\dtau
	\le \frac{\ex^{\nu t'}-1}{\nu}\,c_4\,K_*(\|f\|_{\rm g}^{})^2.
\end{equation}
Bounding the remaining terms in \eqref{q:Wt2} similarly, we find
\begin{equation}\begin{aligned}
   |W^\eps(t)|^2 &\le \ex^{-\nu(t-T_*)}\,|W^\eps(T_*)|^2
		+ \eps\,c_5\,\bigl(K_*^2 + K_*^{3/2}\bigr)\\
	&\le \ex^{-\nu(t-T_*)}\,|W(T_*)|^2
		+ \eps\,c_5\,\bigl(K_*^2 + K_*^{3/2}\bigr).
\end{aligned}\end{equation}
This proves the theorem, with
$K_{\rm g}(\|f\|_{\rm g}^{})^2 = 2\,c_5\,\bigl(K_*^2 + K_*^{3/2}\bigr)$ and
$T_{\rm g}(|\gb W_0|,\|f\|_{\rm g}^{},\eps)=T_*-\log\bigl[\eps\,c_5\,\bigl(K_*+K_*^{1/2}\bigr)\bigr]/\nu$.


\section{Higher-Order Estimates}\label{s:ho}

When $\dy_tf=0$ in the very simple model \eqref{q:1dm}, we can obtain
a better estimate on $x'=x-U$ where $U=\eps f/(\eps\mu+\im)$ than on $x$,
namely that $x'(t)\to0$ as $t\to\infty$;
here $U$ is the (exact, higher-order) {\em slow manifold\/}.
The situation is more complicated when $f$ is time-dependent, or
when $x$ is coupled to a slow variable $y$
with the evolution equations having nonlinear terms.
In this case, it is not generally possible to find $U$ (explicit examples
are known where no such $U$ exists), and thus $x'(t)\not\to0$
as $t\to\infty$ for any $U(y,f;\eps)$.
Nevertheless, it is often possible to find a $U^*$ that gives an
exponentially small bound on $x'(t)$ for large $t$.
We shall do this for the primitive equations.

More concretely, in this section we show that, with reasonable regularity
assumptions on the forcing $f$, the leading-order estimate on the fast variable
$W^\eps$ in the previous section can be sharpened to an exponential-order
estimate on $W^\eps-U^*(W^0,f;\eps)$, where $U^*$ is computed below.
As in \cite{temam-dw:ebal}, we make use of the Gevrey regularity of
the solution and work with a finite-dimensional truncation of the
system, whose description now follows.

Given a fixed $\kappa>0$, we define the low-mode truncation of $W$ by
\begin{equation}\label{q:Pldef}
   W^<(\xb,t) = (\Pl W)(\xb,t)
	:= \sum_{|\kb|<\kappa}^\alpha\, w_\kb^\alpha X_\kb^\alpha \ex^{-\im\omega_\kb^\alpha t/\eps}\ex^{\im\kb\cdot\xb}
\end{equation}
where the sum is taken over $\alpha\in\{0,\pm1\}$
and $\kb\in\Zahl_L$ with $|\kb|<\kappa$.
We also define the high-mode part of $W$ by $W^>:=W-W^<$.
The low- and high-mode parts of the slow and fast variables,
$W^{0<}$, $W^{0>}$, $W^{\eps<}$ and $W^{\eps>}$, are defined in
the obvious manner, i.e.\ $W^{0<}$ with $\alpha=0$ in \eqref{q:Pldef}
and $W^{\eps<}$ with $\alpha\in\{\pm1\}$.
It is clear from \eqref{q:Pldef} and \eqref{q:Hsnorm} that the projection
$\Pl$ is orthogonal in $H^s$, so $\Pl$ commutes with both $A$ and $L$ in
\eqref{q:dW}.
We denote $\Pl B$ by $B^<$.

It follows from the definition that the low-mode part $W^<$ satisfies
a ``reverse Poincar{\'e}'' inequality, i.e.\ for any $s\ge0$,
\begin{equation}\label{q:ipoi}
   |\gb W^<|_{H^s}^{} \le \kappa\,|W^<|_{H^s}^{}\,.
\end{equation}
If $W\in G^\sigma(\Dom)$, the exponential decay of its Fourier coefficients
implies that $W^>$ is exponentially small, that is, for any $s\ge0$,
\begin{equation}\label{q:Wgg}
   |W^>|_{H^s}^{} \le C_s\,\kappa^s\,\ex^{-\sigma\kappa}|W|_{G^\sigma}^{}\,.
\end{equation}
The first inequality evidently also applies to the slow and fast parts
separately, i.e.\ with $W^<$ replaced by $W^{0<}$ or $W^{\eps<}$;
as for \eqref{q:Wgg}, it also holds when $W^>$ on the lhs is replaced
by $W^{0>}$ or $W^{\eps>}$.

\medskip
We recall that the global regularity results of Theorem~0 imply that,
with Gevrey forcing, any solution $W\in H^1(\Dom)$ will be
in $G^\sigma(\Dom)$ after a short time.
As in \cite{temam-dw:ebal} and following \cite{matthies:01}, the central
idea here is to split $W^\eps$ into its low- and high-mode parts.
The high-mode part $W^{\eps>}$ is exponentially small by \eqref{q:Wgg}.
We then compute $U^*(W^{0<},f^<;\eps)$ such that $W^{\eps<}-U^*$
becomes exponentially small after some time.

Following historical precedent in the geophysical literature, it is natural
to present our results in two parts, first locally in time and second globally.
(Here ``local in time'' is used in a sense similar to ``local truncation error''
in numerical analysis, giving a bound on the time derivative of some ``error''.)
The following lemma states that, in a suitable finite-dimensional
space, we can find a ``slow manifold'' $W^{\eps<}=U^*(W^{0<},f^<;\eps)$ on
which the normal velocity of $W^{\eps<}$ is at most exponentially small:

\begin{lemma}\label{t:suba}
Let $s>3/2$ and $\eta>0$ be fixed.
Given $W^0\in H^s(\Dom)$ and $f\in H^s(\Dom)$ with $\dy_tf=0$, there exists
$\eps_{**}(|W^0|_{H^s}^{},|f|_{H^s}^{},\eta)$ such that for $\eps\le\eps_{**}$
one can find $\kappa(\eps)$ and $U^*(W^{0<},f^<;\eps)$
that makes the remainder function
\begin{equation}\label{q:Rsdef}\begin{aligned}
   \Rm^*(W^{0<},f^<;\eps)
	&:= \Pl[(\FD U^*)\,\Gy^*] + \frac1\eps LU^*\\
	&\qquad {}+ B^{\eps<}(W^{0<}+U^*,W^{0<}+U^*) + AU^* - f^{\eps<}
\end{aligned}\end{equation}
exponentially small in $\eps$,
\begin{equation}
   |\Rm^*(W^{0<},f^<;\eps)|_{H^s}^{} \le
	\cnst{r} \bigl[(|W^{0<}|_{H^s}^{} + \eta)^2 + |f|_{H^s}^{}\bigr]\,\exp(-\eta/\eps^{1/4});
\end{equation}
here $\FD U^*$ is the derivative of $U^*$ with respect to $W^{0<}$ and
\begin{equation}
   \Gy^* := -B^{0<}(W^{0<}+U^*,W^{0<}+U^*) - AW^{0<} + f^{0<}.
\end{equation}
\end{lemma}

\bigskip\noindent{\bf Remarks.}

\rmk
The bounds may depend on $s$, $\mu$ and $\Dom$ as well as on $\eta$,
but only the latter is indicated explicitly here and in the proof below.

\rmk
Given $\kappa$ fixed, $U^*$ lives in the same space as $W^{\eps<}$,
that is, $(W^0,U^*)_{L^2}^{}=0$ and $\Pl U^*=U^*$.

\rmk
In the leading-order case of \S\ref{s:o1}, the slow manifold is $U^0=0$
and the local error estimate is incorporated directly into the proof
of Theorem~\ref{t:o1};
we therefore did not put these into a separate lemma.

\rmk\label{r:sm}
Unlike formal constructions in the geophysical literature (see, e.g.,
\cite{ob:97,wbsv:95}), our slow manifold is not defined for all possible
$W^0$ and $\eps$.
Instead, given that $|W^0|_{G^\sigma}^{}\le R$, we can define $U^*$
for all $\eps\le\eps_{**}(R,\sigma)$;
generally, the larger the set of $W^0$ over which $U^*$ is to be defined,
the smaller $\eps$ will have to be.

\rmk
In what follows, we will often write $U^*(W^0,f;\eps)$
for $U^*(\Pl W^0,\Pl f;\eps)$;
this should not cause any confusion.

\bigskip
Using the Lemma and a technique similar to that used to prove
Theorem~\ref{t:o1}, we can bound the ``net forcing'' on $W'=W^{\eps<}-U^*$
by $\Rm^*$.
The dissipation term $AW'$ then ensures that $W'$ eventually decays
to an exponentially small size.
This gives us our global result:

\begin{theorem}\label{t:ho}
Let $W_0\in H^1(\Dom)$ and $\gb f\in G^\sigma(\Dom)$ be given with $\dy_tf=0$.
Then there exist $\eps_*(f;\sigma)$ and $T_*(|\gb W_0|,|\gb f|_{G^\sigma}^{})$
such that for $\eps\le\eps_*$ and $t\ge T_*$, we can approximate the fast
variable $W^\eps(t)$ by a function $U^*(W^0(t),f;\eps)$ of the slow variable
$W^0(t)$ up to an exponential accuracy,
\begin{equation}
   |W^\eps(t)-U^*(W^0(t),f;\eps)|_{L^2}^{}
	\le K_*(|\gb f|_{G^\sigma}^{},\sigma)\,\exp(-\sigma/\eps^{1/4}).
\end{equation}
\end{theorem}

\noindent
As in Theorem~\ref{t:o1}, here $K_*$ is a continuous increasing function
of its arguments; $W(t)=W^0(t)+W^\eps(t)$ is the solution of \eqref{q:dW}
with initial condition $W(0)=W_0$.
As before, the bounds depend on $\mu$ and $\Dom$, but these are not indicated
explicitly.

\bigskip\noindent{\bf Remarks.}

\rmk
With very minor changes in the proof of Theorem~\ref{t:ho} below,
one could also show that, if $f\in H^{n+1}$ and $\dy_tf=0$,
then $|W^\eps(t)-U^n(W^0(t),f;\eps)|_{L^2}^{}$ is bounded as $\eps^{n/4}$
for sufficiently large $n$ and possibly something better for smaller $n$.

\rmk
Recalling remark~\ref{r:sm} above, our slow manifold is only defined for
$\eps$ sufficiently small for a given $|W^{0<}|$ (or equivalently, for
$|W^{0<}|$ sufficiently small for a given $\eps$).
The results of Theorem~0 tell us that $W(t)$ will be inside a ball in
$G^\sigma(\Dom)$ after a sufficiently large $t$;
we use (twice) the radius of this absorbing ball to fix the restriction
on $\eps$.
Thus our approach sheds no light on the analogous problem in
the inviscid case, which has no absorbing set.

\rmk
As proved in \cite{ju:07,kobelkov:06,petcu:3dpe},
assuming sufficiently smooth forcing, the primitive equations admit
a finite-dimensional global attactor.
Theorem~\ref{t:ho} states that, for $\eps\le\eps_*(|f|_{G^\sigma}^{})$,
the solution will enter, and remain in, an exponentially thin neighbourhood
of $U^*(W^{0<},f^<;\eps)$ in $L^2(\Dom)$ after some time.
It follows that the global attractor must then be contained in this
exponentially thin neighbourhood as well.

\rmk
The dynamics on this attractor is generally thought to be chaotic
\cite{temam:iddsmp}.
Thus our present results do not qualitatively affect the finite-time
predictability estimate of \cite{temam-dw:ebal}.

\rmk
When $\dy_tf\ne0$, the slaving relation $U^*$ would have a non-local
dependence on $t$.
Quasi-periodic forcing, however, can be handled by introducing
an auxiliary variable $\boldsymbol{\theta}=(\theta_1,\cdots,\theta_n)$,
where $n$ is the number of independent frequencies of $f$.
The slaving relation $U^*$ would then depend on $\boldsymbol{\theta}$ as well
as on $W^{0<}$.

\rmk
Bounds of this type are only available for the fast variable $W^\eps$;
no special bounds exist for the slow variable $W^0$ except in special cases,
such as when the forcing $f$ is completely fast, $(W^0,f)_{L^2}^{}=0$.


\bigskip
We next present the proofs of Lemma~\ref{t:suba} and Theorem~\ref{t:ho}.
The first one follows closely that in \cite{temam-dw:ebal} which used
a slightly different notation;
we redo it here for notational coherence and since some estimates in
it are needed in the proof of Theorem~\ref{t:ho}.
As before, we write $(\cdot,\cdot)\equiv(\cdot,\cdot)_{L^2}$ and
$|\cdot|\equiv|\cdot|_{L^2}^{}$ when there is no ambiguity.

\subsection{Proof of Lemma~\ref{t:suba}}
As usual, we use $c$ to denote a generic constant which may not be the same
each time it appears.
Constants may depend on $s$ and the domain $\Dom$ (and also on $\mu$ for
non-generic ones), but dependence on $\eta$
is indicated explicitly.
Since $s>3/2$, $H^s(\Dom)$ is a Banach algebra, so if $u$ and $v\in H^s$,
\begin{equation}
   |uv|_s^{} \le c\,|u|_s^{}|v|_s^{}
\end{equation}
where here and henceforth $|\cdot|_s^{} := |\cdot|_{H^s}^{}\,$.
Let us take $\eps\le1$ and $\kappa$ as given for now;
restrictions on $\eps$ will be stated as we go along and $\kappa$ will
be fixed in \eqref{q:deltakappa} below.

We construct the function $U^*$ iteratively as follows.
First, let
\begin{equation}\label{q:U1}
   \frac1\eps LU^1 = - B^{\eps<}(W^{0<},W^{0<}) + f^{\eps<}\,,
\end{equation}
where $U^1\in\textrm{range}\,L$ for uniqueness;
similarly, $U^n\in\textrm{range}\,L$ in what follows.
For $n=1,2,\cdots$, let
\begin{equation}\label{q:Unp1}
   \frac1\eps LU^{n+1} = -\Pl\bigl[(\FD U^n)\Gy^n\bigr]
	- B^{\eps<}(W^{0<}+U^n,W^{0<}+U^n)
	- AU^n + f^{\eps<},
\end{equation}
where $\FD U^n$ is the Fr{\'e}chet derivative of $U^n$ with respect to $W^{0<}$
(regarded as living in an appropriate Hilbert space) and
\begin{equation}
   \Gy^n := -B^{0<}(W^{0<}+U^n,W^{0<}+U^n) - AW^{0<} + f^{0<}.
\end{equation}
We note that the right-hand sides of \eqref{q:U1} and \eqref{q:Unp1} do
not lie in
$\mathrm{ker}\,L$, so $U^1$ and $U^{n+1}$ are well defined.
Moreover, $U^n$ lives in the same space as $W^{\eps<}$, that is,
$U^n\in\Pl\textrm{range}\,L$;
in other words, $(W^0,U^n)=0$ and $\Pl U^n=U^n$.

For $\eta>0$, let $D_\eta(W^{0<})$ be the complex $\eta$-neighbourhood of
$W^{0<}$ in $\Pl H^s(\Dom)$.
With $W^{0<}$ defined by \eqref{q:Pldef}, this is
\begin{equation}\label{q:Deta}\begin{aligned}
   D_\eta(W^0) = \biggl\{ &\hat W^0 : \hat W^0(\xb,t)
	= \sum_{|\kb|<\kappa}\,\hat w_\kb^0 X_\kb^0\ex^{\im\kb\cdot\xb} \quad\textrm{with }\\
	&\quad\hat w_{(k_1,k_2,k_3)}^0 = \hat w_{(k_1,k_2,-k_3)}^0 \textrm{ and }
	\sum_{|\kb|<\kappa}\,|\kb|^{2s}\,|\hat w_\kb^0-w_\kb^0|^2 < \eta^2 \biggr\}.
\end{aligned}\end{equation}
Since $W^0(\xb,t)$ and $X_\kb^0$ are real, $w_\kb^0$ must satisfy
(\ref{q:w0}a), but $\hat w_\kb^0$ in \eqref{q:Deta} need not
satisfy this condition although it must satisfy (\ref{q:w0}b).
We can thus regard $D_\eta(W^{0<}) \subset \{ (w_\kb^{}) : 0 < |\kb| < \kappa
\textrm{ and } w_{(k_1,k_2,-k_3)}=w_{k_1,k_2,k_3)} \} \cong \Comp^m$
for some $m$.
Let $\delta>0$ be given; it will be fixed below in \eqref{q:deltakappa}.
For any function $g$ of $W^{0<}$, let
\begin{equation}
   |g(W^{0<})|_{s;n}^{} := \sup_{W\in D_{\eta-n\delta}(W^{0<})}\,|g(W)|_s^{}\,;
\end{equation}
this expression is meaningful when $D_{\eta-n\delta}(W^{0<})$ is non-empty,
that is, for $n\in\{0,\cdots,\lfloor\eta/\delta\rfloor =: n_*\}$.
For future reference, we note that
\begin{equation}
   |W^{0<}|_{s;0}^{} \le |W^{0<}|_s^{} + \eta.
\end{equation}

Our first step is to obtain by induction a couple of uniform bounds
\eqref{q:bdUn}--\eqref{q:bdUW}, valid for $n\in\{1,\cdots,n_*\}$,
which will be useful later.
First, for $U^1$, we have
\begin{equation}
   \frac1\eps|LU^1|_{s;1}^{} \le |B^{\eps<}(W^{0<},W^{0<})|_{s;1}^{} + |f^{\eps<}|_s^{}
\end{equation}
which, using the estimate $|B(W,W)|_s^{}\le c\,|\gb W|_s^2$ and \eqref{q:ipoi},
implies
\begin{equation}\label{q:bdU1}
   |U^1|_{s;1}^{} \le \eps\,\cnst0\,\bigl(\kappa^2|W^{0<}|_{s;1}^2 + |f^{\eps<}|_{s}^{}\bigr).
\end{equation}

Next, we derive an iterative estimate for $|U^n|_{s;n}^{}$.
Using the fact that $|\cdot|_{s;m}^{} \le |\cdot|_{s;n}^{}$ whenever $m\ge n$,
we have for $n=1,2,\cdots$,
\begin{equation}\begin{aligned}
   \frac1\eps\,|U^{n+1}|_{s;n+1}^{} \le |(\FD U^n)\Gy^n|_{s;n+1}^{}
	&+ |B^{\eps<}(W^{0<}+U^n,W^{0<}+U^n)|_{s;n}^{}\\
	&+ \mu\kappa^2\,|W^{0<}|_{s;n}^{} + |f^{\eps<}|_s^{}\,.
\end{aligned}\end{equation}
The first term on the right-hand side can be bounded by a technique
based on Cauchy's integral formula:
Let $D_\eta(z_0^{})\subset\Comp$ be the complex $\eta$-neighbourhood
of $z_0^{}$.
For $\vfi: D_\eta(z_0^{})\to\Comp$ analytic and $\delta\in(0,\eta)$, we can
bound $|\vfi'|$ in $D_{\eta-\delta}(z_0^{})$ by $|\vfi|$ in $D_\eta(z_0^{})$ as
\begin{equation}\label{q:cauchy}
   |\vfi'\cdot z|_{D_{\eta-\delta}(z_0^{})}
	\le \frac1\delta |\vfi|_{D_\eta(z_0^{})}^{}|z|_{\Comp}^{}\,.
\end{equation}
Now by \eqref{q:U1} $U^1$ is an analytic function of the finite-dimensional
variable $W^{0<}$, so assuming that $U^n$ is analytic in $W^{0<}$ we can
regard the Fr{\'e}chet derivative $\FD U^n$ as an ordinary derivative.
Taking for $\vfi'$ in \eqref{q:cauchy} the derivative of $U^n$ in the
direction $\Gy^n$ (i.e.\ working on the complex plane containing $0$
and $\Gy^n$), we have
\begin{equation}
   |(\FD U^n)\Gy^n|_{s;n+1}^{} \le \frac1\delta\, |U^n|_{s;n}|\Gy^n|_{s;n}^{}\,.
\end{equation}
Using the estimate
\begin{equation}
   |B^{\eps<}(W^{0<}+U^n,W^{0<}+U^n)|_{s;n}^{}
	\le c\,|\gb(W^{0<}+U^n)|_{s;n}^2
	\le c\,\kappa^2|W^{0<}+U^n|_{s;n}^2
\end{equation}
we have
\begin{equation}\label{q:bdUn1}\begin{aligned}
   |U^{n+1}|_{s;n+1}^{}
	&\le \frac{\eps c}\delta\,|U^n|_{s;n}^{} \bigl( c\,\kappa^2\,|W^{0<}+U^n|_{s;n}^2 + \mu\kappa^2\,|W^{0<}|_{s;n}^{} + |f^{0<}|_s^{} \bigr)\\
	&\qquad+ \eps\kappa^2\,c\,|W^{0<}+U^n|_{s;n}^2 + \mu\eps\kappa^2\,|U^n|_{s;n}
	+ \eps\,|f^{\eps<}|_s^{}\,.
\end{aligned}\end{equation}

To complete the inductive step, let us now set
\begin{equation}\label{q:deltakappa}
   \delta = \eps^{1/4}
   \qquad\textrm{and}\qquad
   \kappa = \eps^{-1/4}.
\end{equation}
With this, we have from \eqref{q:bdUn1}
\begin{equation}\label{q:bdUn2}\begin{aligned}
      |U^{n+1}|_{s;n+1}^{}
	&\le \eps^{1/4}\,\cnst1\,|U^n|_{s;n}^{} \bigl( |W^{0<}+U^n|_{s;n}^2 + \mu\,|W^{0<}|_{s;n}^{} + \eps^{1/2}\,|f^{0<}|_s^{} \bigr)\\
	&\qquad+ \eps^{1/2}\,\cnst2\,\bigl(|W^{0<}+U^n|_{s;n}^2 + \mu\,|U^n|_{s;n}
	+ \eps^{1/2}\,|f^{\eps<}|_s^{}\bigr).
\end{aligned}\end{equation}
We require $\eps$ to be such that
\begin{equation}\label{q:eps1}
   \eps^{1/4}\,(\cnst0+\cnst1+\cnst2)\,
	\bigl(|W^{0<}|_{s;0}^2 + \mu\,|W^{0<}|_{s;0}^{} + |f|_s^{} \bigr)
	\le \sfrac14\min\{ 1, |W^{0<}|_s^{} \}
\end{equation}
and claim that with this we have
\begin{equation}\label{q:bdUn}
   |U^n|_{s;n}^{} \le \eps^{1/4}\,\cU\,
	\bigl(|W^{0<}|_{s;0}^2 + \mu\,|W^{0<}|_{s;0}^{} + |f^<|_s^{}\bigr)
\end{equation}
with $\cU=4\,(\cnst0+\cnst1+\cnst2)$.
Now since $\eps\le1$, \eqref{q:bdU1} implies that it holds for $n=1$,
so let us suppose that it holds for $m=0,\cdots,n$ for some $n<n_*$.
Now \eqref{q:eps1} and \eqref{q:bdUn} imply that
\begin{equation}\label{q:bdUW}
   |U^m|_{s;m}^{} \le |W^{0<}|_s^{} \le |W^{0<}|_{s;0}^{}
   \qquad\textrm{and}\qquad
   |U^m|_{s;m}^{} \le 1
\end{equation}
for $m=0,\cdots,n$.
Using these in \eqref{q:bdUn2}, we have
\begin{equation}\begin{aligned}
   |U^{n+1}|_{s;n+1}^{}
	&\le 4\,\eps^{1/4}\,\cnst1\,\bigl(|W^{0<}|_{s;0}^2 + \mu\,|W^{0<}|_{s;0}^{} + |f^<|_s^{}\bigr)\,|U^n|_{s;n}^{}\\
	&\hskip30pt {}+ 4\,\eps^{1/2}\,\cnst2\,\bigl(|W^{0<}|_{s;0}^2 + \mu\,|W^{0<}|_{s;0}^{} + |f^<|_s^{}\bigr)\\
	&\le \eps^{1/4}\cU\,\bigl(|W^{0<}|_{s;0}^2 + \mu\,|W^{0<}|_{s;0}^{} + |f^<|_s^{}\bigr).
\end{aligned}\end{equation}
This proves \eqref{q:bdUn} and \eqref{q:bdUW} for $n=0,\cdots,n_*$.

\medskip
We now turn to the remainder
\begin{equation}
   \Rm^0 := B^{\eps<}(W^{0<},W^{0<}) - f^{\eps<}
\end{equation}
and, for $n=1,\cdots$,
\begin{equation}\label{q:Rndef}
   \Rm^n := \Pl[(\FD U^n)\,\Gy^n] + \frac1\eps LU^n
	+ B^{\eps<}(W^{0<}+U^n,W^{0<}+U^n) + AU^n - f^{\eps<}.
\end{equation}
We seek to show that, for $n=0,\cdots,n_*$, it scales as $\ex^{-n}$.
We first note that by construction $\Rm^n\not\in\textrm{ker}\,L$,
so $L^{-1}\Rm^n$ is well-defined.
Taking $U^0=0$, we have
\begin{equation}\label{q:RUU}
   \Rm^n = \frac1\eps L\,(U^n - U^{n+1}).
\end{equation}
We then compute
\begin{equation}\begin{aligned}
   \Rm^{n+1}
	&= \Pl[(\FD U^{n+1})\,\Gy^{n+1}] + \frac1\eps LU^{n+1}\\
	&\qquad {}+ B^{\eps<}(W^{0<}+U^{n+1},W^{0<}+U^{n+1}) + AU^{n+1} - f^{\eps<}\\
	&= \Pl[(\FD U^{n+1})(\Gy^n+\delta\Gy^n)] + \frac1\eps LU^n - \Rm^n\\
	&\qquad {}+ B^{\eps<}(W^{0<}+U^n,W^{0<}+U^n)
	- \eps\,B^{\eps<}(W^{0<}+U^n,L^{-1}\Rm^n)\\
	&\qquad {}- \eps\,B^{\eps<}(L^{-1}\Rm^n,W^{0<}+U^{n+1})
	+ AU^n - \eps\, AL^{-1}\Rm^n - f^{\eps<}\\
	&= \Pl[(\FD U^n)\,\delta\Gy^n]
	- \eps\,L^{-1}\Pl[(\FD\Rm^n)\,\Gy^{n+1}]
        - \eps\,AL^{-1}\Rm^n\\
	&\qquad {}- \eps\,B^{\eps<}(L^{-1}\Rm^n,W^{0<}+U^{n+1})
	- \eps\,B^{\eps<}(W^{0<}+U^n,L^{-1}\Rm^n),
\end{aligned}\end{equation}
where we have used \eqref{q:RUU} and where
\begin{equation}\begin{aligned}
   \delta\Gy^n &:= \Gy^{n+1} - \Gy^n\\
	&= \eps\,B^{0<}(W^{0<}+U^{n+1},L^{-1}\Rm^n)
	+ \eps\,B^{0<}(L^{-1}\Rm^n,W^{0<}+U^n).
\end{aligned}\end{equation}

To obtain a bound on $\Rm^n$, we compute using \eqref{q:bdUW}
\begin{equation}\label{q:Gyn}\begin{aligned}
   |\Gy^n|_{s;n}^{} &\le c\,\bigl( |\gb(W^{0<}+U^n)|_{s;n}^2
	+ \mu\,|\Delta W^{0<}|_{s;n}^{} + |f^{0<}|_s^{}\bigr)\\
	&\le c\,\kappa^2\,\bigl( |W^{0<}|_{s;0}^2 + \mu\,|W^{0<}|_{s;0}^{}
	+ |f|_s^{}\bigr),
\end{aligned}\end{equation}
as well as
\begin{equation}\begin{aligned}
   |\delta\Gy^n|_{s;n+1}^{}
	&\le \eps\,c\,|\gb(W^{0<}+U^{n+1})|_{s;n+1}^{}|\gb L^{-1}\Rm^n|_{s;n+1}^{}\\
	&\qquad {}+ \eps\,c\,|\gb L^{-1}\Rm^n|_{s;n+1}^{}|\gb(W^{0<}+U^n)|_{s;n}^{}\\
	&\le \eps\kappa^2\,c\,|\Rm^n|_{s;n+1}^{}|W^{0<}|_{s;0}^{}\,.
\end{aligned}\end{equation}
(Note that we can only estimate $\delta\Gy^n$ in $D_{\eta-(n+1)\delta}^{}$
and not in $D_{\eta-n\delta}^{}$;
similarly, since the definition of $\Rm^n$ involves $\FD U^n$, it can only
be estimated in $D_{\eta-(n+1)\delta}^{}$.)
We then have
\begin{equation}\begin{aligned}
   \!\!\!\!\!|\Rm^{n+1}|_{s;n+2}^{}
	&\le |\FD U^n|_{s;n+1}^{}|\delta\Gy^n|_{s;n+1}^{}
	+ \eps\,|L^{-1}\FD\Rm^n|_{s;n+2}^{}|\Gy^{n+1}|_{s;n+1}\\
       	&\qquad {}+ \eps\mu\kappa^2\,|\Rm^n|_{s;n+1}^{}
	+ \eps\,|\gb L^{-1}\Rm^n|_{s;n+1}^{}|\gb(W^{0<}+U^{n+1})|_{s;n+1}^{}\\
	&\qquad+ \eps\,|\gb(W^{0<}+U^n)|_{s;n}^{}|\gb L^{-1}\Rm^n|_{s;n+1}^{}\\
	&\le \frac1\delta|U^n|_{s;n}\,\eps\,\kappa^2\,|\Rm^n|_{s;n+1}|W^{0<}|_{s;0}
	+ c\,\frac\eps\delta\,|\Rm^n|_{s;n+1}^{}|\Gy^{n+1}|_{s;n+1}^{}\\
	&\qquad {}+ 4\,\eps\kappa^2\,|\Rm^n|_{s;n+1}^{}|W^{0<}|_{s;0}^{}
	+ \eps\mu\kappa^2\,|\Rm^n|_{s;n+1}^{}\\
	&\le \eps^{1/4}\,|\Rm^n|_{s;n+1}^{}\,\cnst{e}(|W^{0<}|_{s;0}^2 + \mu\,|W^{0<}|_{s;0}^{} + |f^<|_s^{} + \mu)
\end{aligned}\end{equation}
where for the last inequality we have assumed that
\begin{equation}\label{q:eps1a}
   \eps^{1/4} \le \min\{ \mu/|W^{0<}|_{s;0}^{}, \mu\,\cU/4 \}.
\end{equation}
If we require $\eps$ to satisfy, in addition to $\eps\le1$, \eqref{q:eps1} and
\eqref{q:eps1a},
\begin{equation}\label{q:eps2}
   \eps^{1/4}\,\cnst{e}(|W^{0<}|_{s;0}^2 + \mu\,|W^{0<}|_{s;0}^{} + |f^<|_s^{} + \mu)
	\le \frac1\ex\,,
\end{equation}
we have, for $n=0,1,\cdots,n_*-1$,
\begin{equation}
   |\Rm^{n+1}|_{s;n+2}^{} \le \frac1\ex\,|\Rm^n|_{s;n+1}^{}\,.
\end{equation}
Along with the estimate
\begin{equation}
   |\Rm^0|_{s;1}^{} \le \cnst{r}\, (|W^{0<}|_{s;0}^2 + |f^<|_s^{}),
\end{equation}
taking $n=n_*-1$ leads us to
\begin{equation}\label{q:bdRs}\begin{aligned}
   |\Rm^{n_*-1}|_{H^s}^{} \le |\Rm^{n_*-1}|_{s;n_*}^{}
	&\le \cnst{r}\, (|W^{0<}|_{s;0}^2 + |f^<|_s^{})\,\exp(-\eta/\eps^{1/4})\\
	&\le \cnst{r}\, [(|W^{0<}|_s^{} + \eta)^2 + |f^<|_s^{}]\,\exp(-\eta/\eps^{1/4}).
\end{aligned}\end{equation}
The lemma follows by setting $U^*=U^{n_*-1}$ and taking as $\eps_{**}$
the largest value that satisfies $\eps\le1$, \eqref{q:eps1}, \eqref{q:eps1a}
and \eqref{q:eps2}.

For use later in the proof of Theorem~\ref{t:ho}, we also bound
\begin{equation}\label{q:bdRR}\begin{aligned}
   \bigl|\gb(1-\Pl)&[(\FD U^*)\Gy^*]\bigr|_{L^2}^{}
	\le c\,\ex^{-\sigma\kappa}|(\FD U^*)\Gy^*|_{2,n_*}\\
	&\le c\,\ex^{-\sigma\kappa}\,\frac1\delta\,|U^*|_{2,n_*-1}^{}|\Gy^*|_{2,n_*-1}\\
	&\le c\,\ex^{-\sigma\kappa}\,\kappa^2\,(|W^{0<}|_{2;0}^2 + \mu\,|W^{0<}|_{2;0}^{} + |f|_2^{})^2
\end{aligned}\end{equation}
where for the last inequality we have used \eqref{q:bdUn} and
\eqref{q:Gyn} with $n=n_*-1$.


\subsection{Proof of Theorem~\ref{t:ho}}
We follow the conventions of the proofs of Theorem~\ref{t:o1}
and Lemma~\ref{t:suba} on constants.
We will be rather terse in parts of this proof which mirror a development
in the proof of Theorem~\ref{t:o1}.

First, we recall Theorem~0 and consider
$t\ge T:=\max\{T_2,T_\sigma\}$ so that $|\gb^2W(t)|\le K_2$
and $|\gb^2W(t)|_{G^\sigma}^{}\le M_\sigma$.
We use Lemma~\ref{t:suba} with $s=2$ and, collecting the constraints
on $\eps$ there, require that
\begin{equation}\label{q:eps1*}\begin{aligned}
   &\eps^{1/4}\,\cU\,\bigl( (K_2 + \eta)^2 + \mu\,(K_2 + \eta) + |f|_2^{}\bigr) \le
	\sfrac14 \min\{ 1, K_2 \},\\
   &\eps^{1/4} \le \min\{ \mu/(K_2 + \eta), \mu\,\cU/4, 1 \},\\
   &\eps^{1/4}\, \cnst{e}\,\bigl( (K_2 + \eta)^2 + \mu\,(K_2 + \eta) + \mu + |f|_2^{}\bigr)
	\le \frac1\ex,
\end{aligned}\end{equation}
where $\cnst{e}$ is that in \eqref{q:eps2}.
(We note that all these constraints are convex in $K_2$,
so they do not cause problems when $|W^{0<}|<K_2$.)
Further constraints on $\eps$ will be imposed below.
We note the bound \eqref{q:bdRs} and
\begin{equation}\label{q:bdUs}
   |U^*|_{H^2}^{} \le |U^*|_{2;n_*}^{}
	\le \eps^{1/4}\,\cU\,\bigl( (K_2 + \eta)^2 + \mu\,(K_2 + \eta) + |f|_2^{}\bigr)
\end{equation}
which follows from \eqref{q:bdUW}.

We fix $\kappa=\eps^{-1/4}$ as in \eqref{q:deltakappa}
and consider the equation of motion for the low modes $W^<$,
\begin{equation}\label{q:ddtWl}\begin{aligned}
   \dy_tW^< + \frac1\eps LW^< + B^<(W^<,W^<) &+ AW^< - f^<\\
	&= - B^<(W^>,W) - B^<(W^<,W^>)\\
	&=: \hat{\mathcal{H}}\,.
\end{aligned}\end{equation}
Writing
\begin{equation}
  W^{\eps<} = U^*(W^{0<},f^<;\eps) + W'\,,
\end{equation}
the equation governing the finite-dimensional variable $W'(t)$ is
\begin{equation}\begin{aligned}
   \dy_tW' + \frac1\eps LW' + B^{\eps<}(W^<,W^<) &+ AW'\\
	&= -\dy_t U^* - \frac1\eps LU^* - AU^* + f^{\eps<} + \hat\Hm.
\end{aligned}\end{equation}
Using \eqref{q:Rsdef}, this can be written as
\begin{equation}\label{q:ddtWp}\begin{aligned}
   \dy_tW' &+ \frac1\eps LW' + B^{\eps<}(W^<,W') + B^{\eps<}(W',W^{0<}+U^*) + AW'\\
	&= -\Rm^* - (1-\Pl)[(\FD U^*)\,\Gy^*] + \hat\Hm\\
	&=: -\Rm^* + \Hm.
\end{aligned}\end{equation}

Multiplying by $W'$ in $L^2(\Dom)$, we find
\begin{equation}\label{q:ddt0Wp}
   \frac12\ddt{\;}|W'|^2 + (W',B^{\eps<}(W',W^{0<}+U^*)) + \mu\,|\gb W'|^2
	= -(W',\Rm^*) + (W',\Hm).
\end{equation}
We now write the nonlinear term as
\begin{equation}\begin{aligned}
   (W',B^{\eps<}(W',W^{0<}+U^*))
	&= (W',B(W',W^{0<}+U^*))\\
	&= (W',B(W',U^*)) + (W',B(W',W^{0<}))\\
	&= (W',B(W',U^*)) - (W^{0<},B(W',W')).
\end{aligned}\end{equation}
Following the proof of Theorem~\ref{t:o1} [cf.~\eqref{q:ddtwfour}],
we rewrite \eqref{q:ddt0Wp} as
\begin{equation}\label{q:ddt1Wp}\begin{aligned}
   \ddt{\;}\bigl(\ex^{\nu t}|W'|^2\bigr) &+ \mu\,\ex^{\nu t}\,|\gb W'|^2
	\le -2\,\ex^{\nu t}\,(W',\Rm^*) + 2\,\ex^{\nu t}\,(W',\Hm)\\
	&{}- 2\,\ex^{\nu t}\,(W',B(W',U^*))
	+ 2\,\ex^{\nu t}\,(W^{0<},B(W',W')).
\end{aligned}\end{equation}

We bound the first two terms on the right-hand side as
\begin{equation}\begin{aligned}
   &2\,|(W',\Rm^*)| \le \frac\mu6\,|\gb W'|^2 + \frac{c}\mu\,|\Rm^*|^2,\\
   &2\,|(W',\Hm)| \le \frac\mu6\,|\gb W'|^2 + \frac{c}\mu\,|\Hm|^2.
\end{aligned}\end{equation}
As for the third term in \eqref{q:ddt1Wp}, we bound it as
\begin{equation}\begin{aligned}
   2\,|(W'&,B(W',U^*))| \le c\,|W'|_{L^6}|\gb W'|_{L^2}|\gb U^*|_{L^3}\\
	&\le |\gb W'|^2 \>
	\cnst1\,\eps^{1/4}\bigl( (|W^{0<}|_2^{} + \eta)^2 + \mu\,(|W^{0<}|_2^{} + \eta) + |f^<|_2^{}\bigr)\\
\end{aligned}\end{equation}
where we have used \eqref{q:bdUs} in the last step.
We now require $\eps$ to be small enough so that
\begin{equation}\label{q:eps10}
   \eps^{1/4}\cnst1\,\bigl( (K_2 + \eta)^2 + \mu\,K_2^{} + \mu\,\eta + |f|_2^{}\bigr) \le \frac\mu6,
\end{equation}
which implies that, since $|W^{0<}|_2^{}\le K_2$ by hypothesis,
\begin{equation}
   2\,|(W',B(W',U^*))| \le \frac\mu6\,|\gb W'|^2.
\end{equation}
With these estimates, \eqref{q:ddt1Wp} becomes
\begin{equation}\label{q:ddt2Wp}\begin{aligned}
   \ddt{\;}\bigl(\ex^{\nu t}|W'|^2\bigr) + \frac\mu2\,\ex^{\nu t}\,|\gb W'|^2
	&\le \frac{c}\mu\,\ex^{\nu t}\,\bigl( |\Rm^*|^2 + |\Hm|^2 \bigr)\\
	&\qquad+ 2\,\ex^{\nu t}\,(W^{0<},B(W',W')).
\end{aligned}\end{equation}

Integrating this inequality and multiplying by $\ex^{-\nu T}$, we find
\begin{equation}\label{q:WpTt}\begin{aligned}
   \!\!\ex^{\nu t}\,|W'(T&+t)|^2 - |W'(T)|^2
	+ \frac\mu2 \int_T^{T+t} \ex^{\nu(\tau-T)}|\gb W'|^2 \dtau\\
	&\le \int_T^{T+t} \ex^{\nu(\tau-T)}\,\Bigl\{
		\frac{c}\mu \bigl(|\Rm^*|^2 + |\Hm|^2\bigr)
		+ 2\,(W^{0<},B(W',W')) \Bigr\} \dtau.
\end{aligned}\end{equation}
We then integrate the last term by parts as in \eqref{q:Bws},
\begin{equation}\label{q:IBws}\begin{aligned}
   \!\!\!\!\!\int_T^{T+t} &\ex^{\nu(\tau-T)}\,(W^{0<},B(W',W')) \dtau\\
	&= \eps\,\ex^{\nu t}\,(W^{0<},\Bw(W',W'))(T+t)
	- \eps\,(W^{0<},\Bw(W',W'))(T)\\
	&\quad {}- \eps\int_T^{T+t} \ex^{\nu(\tau-T)}\,\bigl\{
		\nu\,(W^{0<},\Bw(W',W'))
		+ (\dy_\tau W^{0<},\Bw(W',W'))\\
		&\hskip169pt {}+ 2\,(W^{0<},\Bw(\dstau W',W')) \bigr\} \dtau.
\end{aligned}\end{equation}
To bound the terms in the integral, we first need to estimate
\begin{equation}\begin{aligned}
   |\gb B^{\eps<}(W',W^{0<}+U^*)|_{L^2}^{}
	&\le \kappa\,|B^{\eps<}(W',W^{0<}+U^*)|_{L^2}^{}\\
	&\le c\,\kappa\,|\gb W'|_{L^2}^{}|\gb(W^{0<}+U^*)|_{L^\infty}^{}\\
	&\le c\,\kappa^2\,|\gb W'|\,|W^{0<}+U^*|_{H^2}\\
	&\le c\,\kappa^2\,|\gb W'|\,|\gb^2W^0|
\end{aligned}\end{equation}
where for the last inequality we have used \eqref{q:bdUW}.
Using this and the bound
\begin{equation}
   |\gb B^{\eps<}(W^<,W')|_{L^2}^{}
	\le c\,\kappa\,|\gb W^<|_{L^\infty}^{}|\gb W'|_{L^2}
	\le c\,\kappa^2\,|\gb^2 W|\,|\gb W'|
\end{equation}
for the term $B^{\eps<}(W^<,W')$ in \eqref{q:ddtWp} gives us
\begin{equation}
   |\gb\dst W'|_{L^2}^{} \le c\,\kappa^2\,|\gb W'|\,|\gb^2W|
	+ \mu\,\kappa^2\,|\gb W'| + |\gb\Rm^*| + |\gb\Hm|.
\end{equation}
The worst term in \eqref{q:IBws} can now be bounded as
\begin{equation}\begin{aligned}
   \eps\,|(W^{0<},\Bw(\dst W',W'))|
	&\le \eps\,c\,|W^{0<}|_{L^\infty}^{}|\gb W'|_{L^2}^{}|\gb\dst W'|_{L^2}^{}\\
	&\le \cnst2\,\eps\,\kappa^2\,K_2^2\,|\gb W'|^2
	+ \cnst3\,\eps\,\kappa^2\,\mu\,K_2\,|\gb W'|^2\\
	&\qquad {}+ \frac\mu{48}\,|\gb W'|^2
	+ \frac{\eps^2c\,K_2^2}\mu\,(|\gb\Rm^*|^2 + |\gb\Hm|^2).
\end{aligned}\end{equation}
If we now require that $\eps$ satisfy
\begin{equation}\label{q:eps11}
   \eps^{1/2}\,\cnst2\,K_2^2 \le \frac\mu{48}
   \qquad\textrm{and}\qquad
   \eps^{1/2}\,\cnst3\,K_2 \le \frac1{48}\,,
\end{equation}
we have
\begin{equation}
   \eps\,|(W^{0<},\Bw(\dst W',W'))|
	\le \frac\mu{16}\,|\gb W'|^2
	+ \frac{\eps^2c\,K_2^2}\mu\,(|\gb\Rm^*|^2 + |\gb\Hm|^2).
\end{equation}
Bounding another term in \eqref{q:IBws} as
\begin{equation}\begin{aligned}
   \eps\,|(\dy_tW^{0<},\Bw(W',W'))|
	&\le \eps\,c\,|\dy_tW^0|_{L^\infty}^{}|\gb W'|_{L^2}^2\\
	&\le \eps\,c\,|\dy_tW^0|_{H^2}^{}|\gb W'|_{L^2}^2\\
	&\le \eps\,\cnst4\,(\kappa K_2^2 + \mu\kappa^2 K_2 + |f|_2^{})\,|\gb W'|^2
\end{aligned}\end{equation}
and requiring that $\eps$ also satisfy
\begin{equation}\label{q:eps12}
   \eps^{1/2}\,\cnst4\,(K_2^2 + \mu K_2 + |f|_2^{}) \le \frac{\mu}{12},
\end{equation}
plus a similar estimate for the first (easiest) term in \eqref{q:IBws},
we can bound the integral on the r.h.s.\ as
\begin{equation}\begin{aligned}
   \!\int_T^{T+t} &\ex^{\nu(\tau-T)}\,\bigl| (W^{0<},B(W',W')) \bigr| \dtau\\
	&\le \frac\mu2 \int_T^{T+t} \ex^{\nu(\tau-T)}\,|\gb W'|^2 \dtau
	+ \frac{\eps\,c}{\mu^2}\,K_2^2\,(\|\gb\Rm^*\|^2 + \|\gb\Hm\|^2)\,(\ex^{\nu t}-1)
\end{aligned}\end{equation}
where $\|\gb\Rm^*\|:=\sup_{|W^0|\le K_2}|\gb\Rm^*(W^0,f;\eps)|$
and similarly for $\|\gb\Hm\|$.
Bounding the limit term in \eqref{q:IBws} as
\begin{equation}
   |(W^{0<},\Bw(W',W'))|
	\le C\,|W^{0<}|_{L^\infty}^{}|\gb W'|_{L^2}^2
	\le c\,K_2\,\kappa^2\,|W'|^2,
\end{equation}
\eqref{q:WpTt} becomes
\begin{equation}\begin{aligned}
   (1-\eps^{1/2}\cnst5\,&K_2)\,|W'(T+t)|^2\\
	&\le \ex^{-\nu t}\,(1 + \eps^{1/2}\cnst5\,K_2)\,|W'(T)|^2
	+ \frac{\eps\,c}{\mu^2}\,K_2^2\,(\|\gb\Rm^*\|^2 + \|\gb\Hm\|^2).
\end{aligned}\end{equation}


To estimate $\|\gb\Hm\|$, we use \eqref{q:ddtWl}, \eqref{q:Wgg}
and \eqref{q:WGsig}, to obtain
\begin{equation}\begin{aligned}
   |\gb B^{\eps<}(W^>,W)|_{L^2}^{} + |\gb B^{\eps<}(W^<,W^>)|_{L^2}^{}
	&\le c\,\kappa\,|\gb W^>|_{L^4}^{}|\gb W|_{L^4}^{}\\
	&\le c\,\kappa\,\ex^{-\sigma\kappa}\,M_\sigma\,K_2\,.
\end{aligned}\end{equation}
Now \eqref{q:bdRR} implies that
\begin{equation}
   \bigl|\gb(1-\Pl)[(\FD U^*)\Gy^*]\bigr|_{L^2}^{}
	\le c\,\ex^{-\sigma\kappa}\,\kappa^2\,\bigl((K_2 + \eta)^2 + \mu\,(K_2 + \eta) + |f|_2^{}\bigr)^2;
\end{equation}
this and the previous estimate give us
\begin{equation}
   \|\gb\Hm\|_{L^2}^{} \le c\,\ex^{-\sigma\kappa}\,\kappa^2\, \bigl[ M_\sigma K_2
	+ \bigl( (K_2 + \eta)^2 + \mu\,(K_2 + \eta) + |f|_2^{}\bigr)^2 \bigr].
\end{equation}
Meanwhile, using \eqref{q:bdRs} we have
\begin{equation}
   \|\gb\Rm^*\|_{L^2}^{} \le c\,\bigl( (K_2 + \eta)^2 + |f|_2^{}\bigr)\,\exp(-\eta/\eps^{1/4}).
\end{equation}
Setting $\eta=\sigma$ and requiring $\eps$ to satisfy, in addition to
\eqref{q:eps1*}, \eqref{q:eps11} and \eqref{q:eps12},
\begin{equation}\label{q:eps13}
   \eps^{1/2}\,\cnst5\,K_2 \le \frac12\,,
\end{equation}
we have
\begin{equation}\begin{aligned}
   |&W'(T+t)|^2 \le 4\,\ex^{-\nu t}\,|W'(T)|^2\\
	&\quad{}+ \frac{c}{\mu^2}\,\bigl[ (K_2 + \sigma)^4 + \mu^2(K_2+\sigma)^2 + |f|_2^2 + |f|_2^4 + M_\sigma^2 K_2^2
	   \bigr]\,\exp(-2\sigma/\eps^{1/4}).
\end{aligned}\end{equation}
Since $|W'(T)|\le c\,K_2$ by Theorem~0, by taking $t$ sufficiently large
we have
\begin{equation}\begin{aligned}
   |W'(T+t)| \le \frac{c}{\mu}
	\bigl[ (K_2+\sigma)^4 + \mu^2(K_2+\sigma)^2 + |f|_2^2 +{} &|f|_2^4 + M_\sigma^2 K_2^2\bigr] \times\\
	&\eps^{1/2}\,\exp(-\sigma/\eps^{1/4}).
\end{aligned}\end{equation}
And since
\begin{equation}
   |W^\eps-U^*|^2 \le |W^{\eps>}|^2 + |W'|^2
	\le c\,M_\sigma^2\exp(-2\sigma/\eps^{1/4}) + |W'|^2,
\end{equation}
the theorem follows by the same argument used to obtain Theorem~\ref{t:o1}.


\appendix

\section{}\label{s:nores}

\noindent{\bf Proof of Lemma~\ref{t:nores}.}
Since $B_{\jb\kb\lb}^{rs0}=0$ when $j_3k_3|\lb|=0$,
we assume that $j_3k_3|\lb|\ne0$ in the rest of this proof.
As before, all wavevectors are understood to live in $\ZL-\{0\}$ and
their third component take values in $\{0,\pm2\pi/L_3,\pm4\pi/L_3,\cdots\}$.

\medskip
We start by noting that an {\em exact\/} resonance is only possible
when $\jb$ and $\kb$ lie on the same ``resonance cone'', that is,
when $|\jb|/|j_3|=|\kb|/|k_3|$, or equivalently,
when $|\jb'|/|j_3|=|\kb'|/|k_3|$.
There are only two cases to consider:

($\textbf{a}'$)~When $\jb'=\kb'=0$, we have
$B_{\jb\kb\lb}^{rs0}=B_{\kb\jb\lb}^{sr0}=0$.

($\textbf{b}'$)~In the generic case $j_3k_3|\jb'|\,|\kb'|\ne0$, direct
computation using the resonance relation $r|\jb|/j_3+s|\kb|/k_3=0$ gives
$B_{\jb\kb\lb}^{rs0}+B_{\kb\jb\lb}^{sr0}=0$.
This result also follows as the special case $\omega_\jb^r+\omega_\kb^s=0$
in \eqref{q:res0a} below.

\medskip
Now we turn to {\em near\/} resonances.
There are several cases to consider, and we start with the generic
(and hardest) one.

\medskip
($\textbf{a}$)~Suppose that $|\jb'|\,|\kb'|\ne0$ with $\lb'\ne0$.
We define $\Omega$ and $\theta$ by
\begin{equation}
   2\Omega := \omega_\jb^r - \omega_\kb^s
   \qquad\textrm{and}\qquad
   2\theta\Omega := \omega_\jb^r + \omega_\kb^s.
\end{equation}
(We note that $\Omega$ and $\theta$ could take either sign.
Our concern is obviously with small $|\theta|$, when when $\omega_\jb^r$ and
$\omega_\kb^s$ are nearly resonant, so we will restrict $\theta$ below.)
Now this implies that
\begin{equation}
   \omega_\jb^r = (1+\theta)\Omega
   \qquad\textrm{and}\qquad
   \omega_\kb^s = (\theta-1)\Omega.
\end{equation}
We first note that
\begin{equation}
   |\jb'|^2/|j_3|^2 = (1+\theta)^2\Omega^2 - 1
   \qquad\textrm{and}\qquad
   |\kb'|^2/|k_3|^2 = (1-\theta)^2\Omega^2 - 1
\end{equation}
and compute
\begin{equation}
   |\jb'|^2\frac{k_3}{j_3} - |\kb'|^2\frac{j_3}{k_3}
	= 4\theta\Omega^2j_3k_3
	= -\frac{4\theta}{1-\theta^2}\,rs\,|\jb|\,|\kb|.
\end{equation}

Direct computation gives us
\begin{equation}\begin{aligned}
   B_{\jb\kb\lb}^{rs0} &+ B_{\kb\jb\lb}^{sr0}\\
	&= \frac{\im|\Dom|\delta_{\jb+\kb-\lb}\,j_3k_3}{2\,|\jb|\,|\jb'|\,|\kb|\,|\kb'|\,|\lb|}
	\bigl[(P+P')(Q+Q') + (-P+P'')(Q+Q'')\bigr]\\
	&= \frac{\im|\Dom|\delta_{\jb+\kb-\lb}\,j_3k_3}{2\,|\jb|\,|\jb'|\,|\kb|\,|\kb'|\,|\lb|}
	\bigl[P(Q'-Q'') + (P'+P'')Q + P'Q' + P''Q''\bigr]
\end{aligned}\end{equation}
where
\begin{equation}\begin{aligned}
   &P := \jb'\wedge\kb'
   &&Q := -\jb'\cdot\kb'\\
   &P' := \im\frac{r|\jb|}{j_3}(\jb'\cdot\kb')
	- \im\frac{r|\jb|}{j_3}\frac{k_3}{j_3}|\jb'|^2 \qquad
   &&Q' := -\im\frac{s|\kb|}{k_3}(\jb'\wedge\kb') + \frac{j_3}{k_3}|\kb'|^2\\
   &P'' := \im\frac{s|\kb|}{k_3}(\jb'\cdot\kb')
	- \im\frac{s|\kb|}{k_3}\frac{j_3}{k_3}|\kb'|^2
   && Q'' := \im\frac{r|\jb|}{j_3}(\jb'\wedge\kb') + \frac{k_3}{j_3}|\jb'|^2.
\end{aligned}\end{equation}
After some computation, we find
\begin{equation}\begin{aligned}
   &P' + P'' = 2\theta\Omega\,\im\,\Bigl(
	\jb'\cdot\kb' + \frac{2rs}{1-\theta^2}|\jb|\,|\kb|
	- \frac{|\jb'|^2}2\frac{k_3}{j_3} - \frac{|\kb'|^2}2\frac{j_3}{k_3}
	\Bigr),\\
   &Q'-Q'' = 2\theta\Omega\,\Bigl(
	\frac{2rs/\Omega}{1-\theta^2}\,|\jb|\,|\kb| - \im(\jb'\wedge\kb') \Bigr),\\
   &P'Q' + P''Q'' = 2\theta\Omega\,\Bigl\{
	\frac{2rs\,|\jb|\,|\kb|}{1-\theta^2}\frac{r|\jb|}{j_3}\frac{s|\kb|}{k_3}
			(\jb'\wedge\kb')
	+ 2\im rs|\jb|\,|\kb|(\jb'\cdot\kb')\\
	&\hbox to110pt{}+ \frac{\im}2(\jb'\cdot\kb')\,\Bigl(|\kb'|^2\frac{j_3}{k_3} + |\jb'|^2\frac{k_3}{j_3}\Bigr)
	- \im\,|\jb'|^2\,|\kb'|^2\Bigr\},
\end{aligned}\end{equation}
from which we obtain
\begin{equation}\begin{aligned}
   B_{\jb\kb\lb}^{rs0} &+ B_{\kb\jb\lb}^{sr0}
	= 2\theta\Omega\,\frac{\im\,|\Dom|\,\delta_{\jb+\kb-\lb}}{2\,|\lb|}
	\Bigl\{ \im\,(\jb'\cdot\kb')\,\frac{|\jb'|^2k_3^2+|\kb'|^2j_3^2}{|\jb|\,|\jb'|\,|\kb|\,|\kb'|}
	- 2\im\frac{|\jb'|\,|\kb'|\,j_3k_3}{|\jb|\,|\kb|}\\
	&\!\!- \frac{2\im\,rs\,\theta^2}{1-\theta^2}\frac{(\jb'\cdot\kb')j_3k_3}{|\jb'|\,|\kb'|}
	+ \frac{2\,rs/\Omega}{1-\theta^2}\frac{\jb'\wedge\kb'}{|\jb'|\,|\kb'|}\,j_3k_3
	+ \frac{2/\Omega}{1-\theta^2}\frac{|\jb|\,|\kb|}{|\jb'|\,|\kb'|}(\jb'\wedge\kb') \Bigr\}.
\end{aligned}\end{equation}
Now if we require that $|\theta|\le\theta_0<1$, we have the bound
\begin{equation}\label{q:res0a}
   \bigl|B_{\jb\kb\lb}^{rs0} + B_{\kb\jb\lb}^{sr0}\bigr|
	\le \frac{|\Dom|}{2}\Bigl(4 + \frac{6}{1-\theta_0^2}\Bigr)
		\frac{|\jb|\,|\kb|}{|\lb|}\,|\omega_\jb^r+\omega_\kb^s|.
\end{equation}

To take care of the case $|\theta|>\theta_0$, we note that in this case
\begin{equation}
   |\omega_\jb^r+\omega_\kb^s|
	\ge \theta_0\,\Bigl(\frac{|\jb|}{|j_3|} + \frac{|\kb|}{|k_3|}\Bigr).
\end{equation}
We note that since $\theta_0<1$ by hypothesis, this inequality holds both
when $\omega_\jb^r\omega_\kb^s<0$ and $\omega_\jb^r\omega_\kb^s>0$.
Using \eqref{q:bdBrs0}, we then find
\begin{equation}
   \bigl|B_{\jb\kb\lb}^{rs0} + B_{\kb\jb\lb}^{sr0}\bigr|
	\le \bigl|B_{\jb\kb\lb}^{rs0}\bigr| + \bigl| B_{\kb\jb\lb}^{sr0}\bigr|
	\le \sqrt5\,|\Dom|\,\Bigl(|\kb'| + |\jb'| + |\kb'|\frac{|j_3|}{|k_3|}
		+ |\jb'|\frac{|k_3|}{|j_3|}\Bigr).
\end{equation}
Putting these together, we find after a short computation,
\begin{equation}\label{q:res0b}
   \bigl|B_{\jb\kb\lb}^{rs0} + B_{\kb\jb\lb}^{sr0}\bigr|
	\le \frac{2\sqrt5\,|\Dom|}{\theta_0}\bigl(|j_3|+|k_3|\bigr)
		|\omega_\jb^r+\omega_\kb^s|.
\end{equation}

\medskip
($\textbf{b}$)~Suppose now that $|\jb'|\,|\kb'|\ne0$ but $\lb'=0$.
We find using $\jb'+\kb'=0$,
\begin{equation}
   B_{\jb\kb\lb}^{rs0} + B_{\kb\jb\lb}^{sr0}
	= \im\,|\Dom|\,\delta_{\jb+\kb-\lb}\,
	\frac{-\im\sgn l_3\,|\jb'|\,|\kb'|}{2\,|\jb|\,|\kb|}
	(j_3+k_3)(\omega_\jb^r+\omega_\kb^s),
\end{equation}
and thus the bound
\begin{equation}\label{q:res1}
   \bigl|B_{\jb\kb\lb}^{rs0} + B_{\kb\jb\lb}^{sr0}\bigr|
   \le \frac{|\Dom|}2\,\bigl(|j_3|+|k_3|\bigr)\,|\omega_\jb^r+\omega_\kb^s|.
\end{equation}

\medskip
($\textbf{c}$)~Finally, we consider the case $\jb'=0$ and $\kb'\ne0$
(which obviously implies the case $\kb'=0$ and $\jb'\ne0$).
After some computation using $\lb'=\kb'$, we find
\begin{equation}
   B_{\jb\kb\lb}^{rs0} + B_{\kb\jb\lb}^{sr0}
	= \frac{\im\,|\Dom|\,\delta_{\jb+\kb-\lb}}{2\,|\lb|\,|\kb|}
	j_3(k_1-\im rk_2)|\kb'|\,\Bigl(sr-\frac{|\kb|}{k_3}\Bigr).
\end{equation}
But since in this case
\begin{equation}
   |\omega_\jb^r - \omega_\kb^s|
	= \bigl|r\sgn j_3 - s|\kb|/k_3\bigr|
	= \bigl|rs - |\kb|/k_3\bigr|,
\end{equation}
we have the bound
\begin{equation}\label{q:res2}
   \bigl|B_{\jb\kb\lb}^{rs0} + B_{\kb\jb\lb}^{sr0}\bigr|
	\le \frac{|\Dom|\,|j_3|\,|\kb'|^2}{\sqrt2\,|\kb|\,|\lb|}
		|\omega_\jb^r+\omega_\kb^s|
	\le \frac{|\Dom|\,|j_3|}{\sqrt2}\,|\omega_\jb^r+\omega_\kb^s|,
\end{equation}
which holds whether or not $l_3=0$.
We recall that there is nothing to do when $\jb'=\kb'=0$
since then $B_{\jb\kb\lb}^{rs0}=B_{\kb\jb\lb}^{sr0}=0$.

The lemma follows upon fixing $\theta_0$ and collecting
\eqref{q:res0a}, \eqref{q:res0b}, \eqref{q:res1} and \eqref{q:res2}.


\nocite{temam-ziane:03}


\end{document}